\documentclass{article}
\usepackage{hyperref}
\usepackage[utf8]{inputenc}
\usepackage[cyr]{aeguill}
\usepackage[all]{xy}
\usepackage[french, english]{babel}
\usepackage[T1]{fontenc}
\usepackage{amscd}  
\usepackage{amsmath}
\usepackage{amsthm} 
\usepackage{amssymb}
\usepackage{titlesec}
\usepackage{enumitem}
\theoremstyle{definition}
\newtheorem{theointro}{Theorem}

\newtheorem{propintro}[theointro]{Proposition}   
\newtheorem{theo}{Theorem}[subsubsection]
\newtheorem{defi}[theo]{Definition}
\newtheorem*{defi*}{Definition}
\newtheorem*{Note*}{Notation}
\newtheorem{defis}[theo]{Definitions}                  
\newtheorem{lemm}[theo]{Lemma}
\newtheorem{fact}[theo]{Fact}
\newtheorem{prop}[theo]{Proposition}
\newtheorem{coro}[theo]{Corollary}
\theoremstyle{remark}
\newtheorem{ex}[theo]{Example}
\newtheorem*{ex*}{Example}
\newtheorem{rema}[theo]{Remark}
\newtheorem{note}[theo]{Notation}

\renewcommand\thesubsubsection{\ifnum\arabic{subsection}=0\relax{\thesection}\else{\ifnum\arabic{subsubsection}=0\relax{\thesubsection}\else{\thesubsection.\arabic{subsubsection}}\fi}\fi}
\titleformat{\section}
     {\center\Large\scshape}{\thesection.}{1em}{}
 \titleformat{\subsection}
     {\center\large\scshape}{\thesubsection.}{1em}{}[]
 \titleformat{\subsubsection}
     {\scshape}{\thesubsubsection.}{1em}{}
\def\p{\partial}
\def\scr{\mathcal}                    
\def\tcirc{\,\tilde\circ\,}                    
\def\btimes{\boxtimes}                    
\def\tlambda{\,\tilde\lambda\,}                    
\def\tLambda{\,\tilde\Lambda\,}                    
\def\trho{\,\tilde\rho\,}                                        
\def\c{\operatorname{c}_p}
\def\as{\mathfrak a}
\def\E{\operatorname{E}}       
\def\ad{\operatorname{ad}}       
\def\Sym{\operatorname{Sym}}       
\def\As{\operatorname{As}}       
\def\Lie{\operatorname{Lie}}       
\def\Ind{\operatorname{Ind}}
\def\Shift{\operatorname{Shift}}

\def\Pois{\operatorname{Pois}}
\def\Com{\operatorname{Com}}

\def\id{\operatorname{id}}
\def\N{\mathbb N}
\def\diag{\shorthandoff{;:!?}
\xymatrix}
\def\0{\mathbf 0}
\def\_{\underline}
\def\t{\mathfrak t}
\def\P{\mathcal P}
\def\Lev{\operatorname{Lev}}
\def\O{\scr O}
\def\Q{\scr Q}
\def\F{\mathbb F}
\def\M{\mathcal M}
\def\ı{\mathcal N}
\def\r{{\_r}}
\def\q{{\_q}}
\def\C{\scr C}
\def\∑{\mathfrak S}
\def\to{\rightarrow}
\def\inj{\hookrightarrow}
\def\surj{\twoheadrightarrow}
\def\val#1{\vert#1\vert}
\def\D{\operatorname{D}}

\def\Der{\operatorname{Der}}
\def\Comp{\operatorname{Comp}}
\author{Sacha Ikonicoff}
\begin{document}
	\title{Divided power algebras and distributive laws}
	\maketitle
  \begin{abstract}
    We study divided power structures over a product of operads with distributive law. We give a systematic method to characterise the divided power algebras over such a product from the structures of divided power algebra coming from each of the the factor operads. We characterise divided power algebras with operadic derivation, as well as divided power $p$-level algebras in characteristic $p$, and divided power Poisson algebras in characteristic 3.
  \end{abstract}
	\tableofcontents
  \section{Introduction}
  Divided power algebras were introduced by Cartan \cite{HC} to study the homology of Eilenberg--MacLane spaces with coefficients in a finite field. In \cite{BF}, Fresse defined a more general notion of divided power algebra over an operad. Recall that, if $\P$ is a symmetric operad in vector spaces, then the classical Schur functor associated to $\P$ is the endofunctor $S(\P,-)$ of the category of vector spaces such that:
  \[
    S(\P,V)=\bigoplus_{n\ge 0}\P(n)\otimes_{\∑_n}V^{\otimes n},
  \]
  where $\∑_n$ denotes the symmetric group on $n$ letters, and this group acts diagonally on $\P(n)$ and on the $n$-fold tensor product of $V$ by permuting factors. This functor is endowed with a monad structure, and the associated algebras are the $\P$-algebras. The following gives another endofunctor of the category of vector space:
  \[
     \Gamma(\P,V)=\bigoplus_{n\ge 0}\left(\P(n)\otimes V^{\otimes n}\right)^{\∑_n}.
  \]
   In \cite{BF}, Fresse shows that, when $\P$ is reduced (that is, when $\P(0)=\0$), then the induced functor $\Gamma(\P,-)$ is endowed with a structure of monad. The resulting algebras are called $\Gamma(\P)$-algebras, divided power $\P$-algebras or divided $\P$-algebras. The divided power algebras over the classical operad of commutative, associative algebras are exactly the divided power algebras studied by Cartan in \cite{HC}. Divided power algebras appear classically in crystalline cohomology in positive characteristics \cite{PiB}, and in deformation theory \cite{ThF}. Recent works on deformation theory involve divided power structures constructed from the operad of Lie and PreLie algebras, see \cite{AC} and \cite{BCN}.

   Given an operad $\P$, the structure of $\Gamma(\P)$-algebras is notoriously difficult to characterise. In \cite{BF}, Fresse gives such a characterisation for the operad of associative algebras, for the operad of Lie algebras, and a partial characterisation for the operad of Poisson algebras. In \cite{AC}, Cesaro gives a characterisation of the divided power PreLie algebras. In \cite{SI}, we characterise divided power $\P$-algebras in terms of monomial operations indexed by elements of $\P$ and compositions of integers. This characterisation can be refined by fixing either the operad $\P$ or the characteristic of the base field. A detailed definition of operads and their (divided power) algebras is given in Section \ref{modop}.

   The next step in this investigation is to understand whether divided power structures on an operad can be broken down into simpler structures when the given operad is itself built from simpler operads. One way in which operads can be assembled to form a more complex operad is through the use of distributive laws. Given two operads $\P$ and $\Q$, a distributive law between $\P$ and $\Q$ is a map $\Lambda:\Q\circ\P\to\P\circ\Q$ endowing $\P\circ\Q$ with an operad structure that respects both the operad structure of $\P$ and $\Q$ (see Definition \ref{defDL} for the precise definition).

   In the present article, we obtain the following result:
   \begin{propintro}[\ref{propdistrgamma}]
    Let $\P,\Q$ be two reduced operads, and $\Lambda$ a distributive law on $\P\circ\Q$. A divided power $\P\circ\Q$-algebra is a vector space endowed with both a divided power $\P$-algebra and a divided power $\Q$-algebra structures, with compatibility relations induced by $\Lambda$.
   \end{propintro}
  This result does not come as a surprise and was hinted at in the work of Fresse \cite{BF}. The real contribution of this article is a systematic method to compute divided power $\P\circ\Q$ structures in terms of monomial operations coming from the divided power $\P$ and $\Q$ structures, which is compatible with the results of \cite{SI}. This method is embodied by the computational result of Lemma \ref{lemmtLambda}. The main goal of this article is to develop applications of this method on various examples, which we obtain in Theorems \ref{theodivshiftintro}, \ref{theodivderintro}, \ref{theolevpdivintro} and \ref{theoPois3intro}.

  A motivating example of operad obtained from a distributive law is given by the operad $\Pois$ of Poisson algebra. Recall that a Poisson algebra is a vector space with both a commutative associative multiplication and a Lie bracket, satisfying a certain distributivity relation. The operad $\Pois$ is obtained by equipping the product $\Com\circ\Lie$, where $\Com$ is the operad of associative commutative algebras and $\Lie$ is the operad of Lie algebras, with a distributive law which reflects the distributivity relations of Poisson algebras (see \cite{MMa}). Given an operad $\P$, we will introduce the examples of the operads $\Shift_\P$ and $\Der_\P$, both obtained by distributive laws, which respectively encode the structure of a $\P$-algebra with a $\P$-endomorphism, and the structure of a $\P$-algebra with $\P$-derivation (see Section \ref{secshiftder} for a more precise definition). We introduce a notion of divided $\P$-derivation (see Definition \ref{defdivder}) which generalises the notion of derivation with power rule of Keigher--Pritchard \cite{KP}, and appears in the original work on divided power algebras of Cartan \cite{HC}. Applying Proposition \ref{propdistrgamma} to the operads $\Shift_\P$ and $\Der_\P$ yields the following characterisations for divided $\Shift_\P$- and $\Der_\P$-algebras:
  \begin{theointro}[\ref{theodivshift}]\label{theodivshiftintro}
    A divided power $\Shift_\P$-algebra is a divided power $\P$-algebra endowed with an endomorphism of divided power $\P$-algebras.
  \end{theointro}
\begin{theointro}[\ref{theodivder}]\label{theodivderintro}
    A divided power $\Der_\P$-algebra is a divided power $\P$-algebra endowed with a divided $\P$-derivation.
  \end{theointro}
  We then study divided power algebras over an operad $\Lev_p$ that can be seen as a suboperad of $\Shift_{\Com}$. The operad $\Lev_p$ was introduced in \cite{SI3} to study unstable modules and algebras over the mod $p$ Steenrod algebra. An algebra over the operad $\Lev_p$ is a vector space $A$ endowed with a completely symmetric $p$-ary operation $\star$ satisfying the interchange law (see Definition \ref{exop} for a more precise definition). We obtain the following:
  \begin{theointro}[\ref{theolevpdiv}]\label{theolevpdivintro}
		In characteristic $p$, a divided power $\Lev_p$-algebra is a $\Lev_p$-algebra $(A,\star)$ endowed with an operation $\varphi_p:A\to A$ satisfying, for all $a,b\in A$, and for all scalars $\lambda$,
		\begin{enumerate}[label=\arabic*)]
			\item $\varphi_p(\lambda a)=\lambda^p\varphi_p(a)$,
			\item $\star(a^{\times p})=0$,
			\item $\varphi_p(a+b)=\varphi_p(b)+\left(\sum_{i=1}^{p-1}\frac{\star(a^{\times i},b^{\times {p-i}})}{i!(p-i)!}\right) +\varphi_p(a)$,
			\item $\varphi_p(\star(a_1,\dots,a_p))=0$.
		\end{enumerate}
	\end{theointro}
  In the last section of this article, we study the divided power structure over the operad of Poisson algebras. In this case, Fresse obtained a general characterisation in terms of divided power algebras and divided power Lie-algebras with compatibility relations \cite{BF}. However, one of these relations is only explicit in characteristic 2. In the general case, it involves products of $\Lie$-polynomials $\Gamma_{i,j}$ that appeared in the work of Cohen on the homology of $C_{n+1}$-spaces \cite{FrC}. Using the method developed in this article, we run the explicit computation for these polynomials in characteristic 3 and we obtain:
  \begin{theointro}[\ref{theoPois3}]\label{theoPois3intro}
  In characteristic $p=3$, a $\Gamma(\Pois)$-algebra is a vector space $A$ endowed with a commutative associative multiplication $*$, a Lie bracket $[-;-]$, a divided cube operation $\gamma_3$, a Frobenius $F$, such that $(A,*,\gamma_3)$ is a divided power algebra in the classical sense, $(A,[-;-],F)$ is a $3$-restricted Lie algebra, $(A,*,[-;-])$ is a Poisson algebra, satisfying the compatibility relations which include the following, for all $a,b\in A$:
    \begin{enumerate}[label=(P\arabic*)]
      \item $[\gamma_3(a);b]=\frac{a^{*2}}{2}*[a;b]$,
      \item $F(a*b)=\left[[b;a];a\right]*a*b^{*2}+\left[[a;b];b\right]*a^{*2}*b-a*b*[a;b]^{*2}$,
      \item $F(\gamma_3(a))=0$.
    \end{enumerate}
  \end{theointro}
  We expect that a similar characterisation for divided power Poisson algebra in characteristic $p>3$ may be obtained by the same method, using techniques of computer algebra such as those used by Bremner--Markl \cite{BrMa} and Bremner--Dotsenko \cite{BD}.
\paragraph*{Outline of the paper:}
Sections \ref{modop} to \ref{secchardiv} contain the theoretical background for the study of divided power algebras and the notion of distributive laws on product of operads. Section \ref{secdistrlawdiv} makes the link between these two notions and contain the results (Proposition \ref{propdistrgamma} and Lemma \ref{lemmtLambda}), which allow a systematic study of the divided power structures over a product of operad with distributive laws. Sections \ref{secdivshiftder} to \ref{secPoisdiv} are devoted to examples of application of these results and methods. More precisely:

In Section \ref{modop}, we recall the definitions of operads and their algebras, and divided power algebras over an operad. In Section \ref{secpartcomp}, we build some operations on partitions and compositions of integers that are used in \cite{SI} to obtain the characterisation of divided power algebras over an operad. In Section \ref{secDL}, we review the definition of a distributive law on a product of monads and on a product of operads. In Section \ref{secshiftder}, we use distributive laws to build two operads, $\Shift_{\P}$ and $\Der_{\P}$, which encode respectively the structure of $\P$-algebra with $\P$-algebra endomorphism, and the structure of $\P$-algebra with $\P$-derivation. In Section \ref{secchardiv}, we recall the characterisation of divided $\P$-algebras obtained in \cite{SI} in terms of monomial operations, and we give the first applications of this result to classical operads. In Section \ref{secdistrlawdiv}, we obtain our first result, Proposition \ref{propdistrgamma}. We explain in this section how to express distributive laws over divided power algebra monads in terms of monomial divided power operations. In Section \ref{secdivshiftder}, we apply Proposition \ref{propdistrgamma} to obtain a characterisation of divided power algebras over our operads $\Shift_{\P}$ and $\Der_{\P}$. These characterisations are contained in Theorems \ref{theodivshift} and \ref{theodivder}, respectively. In Section \ref{secdivplev}, we obtain a characterisation of divided $\Lev_p$-algebras, using the fact that $\Lev_p$ injects as an operad in $\Shift_{\Com}$. This is Theorem \ref{theolevpdiv}. Finally, in Section \ref{secPoisdiv}, we recall the characterisations of divided power Poisson algebras due to Fresse \cite{BF}, and we give, in Theorem \ref{theoPois3}, an explicit characterisation of these algebras in characteristic 3.

\paragraph*{Acknowledgements}
The author would like to thank Muriel Livernet and Kristine Bauer Benoît Fresse, and Vladimir Dotsenko, for useful discussions.
\begin{Note*}
  Throughout this article, we work over a base field denoted by $\F$, and $\otimes$ is the tensor product over this field.
\end{Note*}
	\section{Symmetric sequences, operads, algebras, divided power algebras}\label{modop}

	The goal of is section is to set our notation regarding symmetric sequences, operads and their algebras. We introduce the notion of divided power algebra over an operad, due to Fresse \cite{BF}. We give definitions for the functor $\Gamma$, as well as the natural transformation called the trace map, which links the classical Schur functor to the functor $\Gamma$. Both these functors and the trace map play a key role throughout the article. At the end of this section, we give the classical examples of the operads $\As$, $\Com$, $\Lie$, as well as two operads $\D$ and $\Lev_p$ that we introduced in \cite{SI3}. The operad $\D$ is used to construct our operads $\Shift_{\P}$, and $\Der_{\P}$ in Section \ref{secshiftder}. In Section \ref{secdivplev}, we characterise the structure of divided power $\Lev_p$-algebras.

  We assume that the reader is familiar with the notion of symmetric sequences of vector spaces as defined in \cite{LV}. We recall the tensor product and composition product of symmetric sequences, and the more exotic product $\tcirc$ used to define divided power algebras:




\begin{defi}[\cite{LV}, \cite{BF}]
For $\M,\ı$ two symmetric sequences of vector spaces, the tensor product $\M\btimes\ı$ is given by:
  \[
    (\M\btimes\ı)(n):=\bigoplus_{i+j=n}\Ind_{\∑_i\times\∑_j}^{\∑_n}\M(i)\otimes\ı(j).
  \]
  For $\M,\ı$ two symmetric sequences, the symmetric sequences $\M\circ \ı$ and $\M \tcirc\ı$ are given by:
    \[
      \M\circ\ı:=\bigoplus_{k\in\N}\left(\M(k)\otimes \ı^{\btimes\,k}\right)_{\∑_k},
    \]
   and:
   \[
    \M\tcirc\ı:=\bigoplus_{k\in\N}\left(\M(k)\otimes \ı^{\btimes\,k}\right)^{\∑_k},
   \]
    where $\∑_k$ acts diagonally on $\M(k)\otimes \ı^{\btimes\,k}$.

    We denote by $(\mu;\nu_1,\dots,\nu_n)\in \M\circ\ı$ the element $[\mu\otimes(\nu_1\otimes\dots\otimes\nu_n)]$, where $\mu\in\M(n)$ and $\nu_1,\dots,\nu_n\in\ı$.

    There are two symmetric monoidal category structures on $\Sym$, given by the operations $\circ$ and $\tcirc$. The unit for both of these products, denoted by $\scr U$, is the symmetric sequence which is $\F$ in arity 1 and $0$ in any other arity.

    A symmetric sequence $\M$ is said to be reduced (or connected) if $\M(0)=0$.
  \end{defi}
  \begin{rema}
    The natural associator isomorphism $\M\circ(\ı\circ \scr L)\to(\M\circ\ı)\circ\scr L$ maps the element 
    \[
      (\mu;(\nu_1,\lambda_{1,1},\dots,\lambda_{1,n_1})\dots, (\nu_m;\lambda_{m,1},\dots,\lambda_{m,n_m}))
    \] 
    to the element
    \[
      ((\mu;\nu_1,\dots,\nu_m);\lambda_{1,1},\dots,\lambda_{1,n_1},\dots,\lambda_{n,m_n}),
    \]
    for all $\mu\in\M(m)$, $\nu_i\in\ı(n_i)$ for all $i\in\{1,\dots,m\}$, and $\lambda_{i,j}\in\scr L$ for all $i\in\{1,\dots,m\}$ $j\in\{1,\dots,n_i\}$. We will allow the following notation for both of these elements:
    \[
      (\mu;\nu_1,\dots,\nu_m;\lambda_{1,1},\dots,\lambda_{1,n_1},\dots,\lambda_{n,m_n})\in\M\circ\ı\circ\scr L.
    \]
    The natural associator isomorphism $\M\tcirc(\ı\tcirc \scr L)\to(\M\tcirc\ı)\tcirc\scr L$, which will be here denoted by $\as_{\M,\ı,\scr L}$, is less straightforward. We will give an explicit formula for this isomorphism in Lemma \ref{lemass}.
  \end{rema}
  \begin{defi}[Schur functor, divided Schur functor, trace map]
    Let $\M$ be a symmetric sequence. The Schur functor associated to $\M$ is the endofunctor of the category of $\F$-vector spaces, mapping $V$ to
  \[
    S(\M,V)=\bigoplus_{n\in\N}(\M(n)\otimes V^{\otimes n})_{\∑_n}=\M\circ V,
  \]
  where $V$ is identified to a symmetric sequence concentrated in arity 0. This endofunctor is denoted $S(\M)$.

  Similarly, for all symmetric sequences $\M$, there is an endofunctor of the category of $\F$-vector spaces, which maps $V$ to
   \[
    \Gamma(\M,V)=\bigoplus_{n\in\N}(\M(n)\otimes V^{\otimes n})^{\∑_n}=\M\tcirc V.
   \]
    This endofunctor is denoted $\Gamma(\M)$.

    For all symmetric sequences $\M$ and $\ı$, the norm map or the trace map is the morphism of symmetric sequences, natural in $\M$ and in $\ı$:
  \[
    Tr_{\M,\ı}:\M\circ \ı\to\M\tcirc\ı,
  \]
  This induces a natural transformation, also called the trace, $Tr_\M:S(\M)\to\Gamma(\M)$.
  \end{defi}
  \begin{lemm}[\cite{BF}]\label{lemmTriso}
  If $\ı$ is reduced, then $Tr_{\M,\ı}$ is an isomorphism.
  \end{lemm}
  \begin{defi}[Operad, algebras over an operad, divided power algebras over an operad]\label{defiop}
    An operad $\P$ is a monoid in the monoidal category $(\Sym,\circ,\scr U)$. In other words, $\P$ is a symmetric sequence endowed with morphisms $\eta_\P:\scr U\to\P$ and $\mu_\P:\P\circ\P\to\P$ satisfying unitality, equivariance and associativity axioms. We denote by $1_\P\in\P(1)$ the image of $1\in\F$ by $\eta_\P$, and by $\mu(\nu_1,\dots,\nu_n)$ the image of $(\mu;\nu_1,\dots,\nu_n)$ by $\mu_\P$.

  If $\P$ is an operad, the endofunctor $S(\P)$ is a monad, and we define $\P$-algebras to be the algebras over the monad $S(\P)$. If $A$ is a $\P$-algebra, $\mu\in \P(n)$, $a_1,\dots,a_n\in A$, we denote by $\mu(a_1,\dots,a_n)\in A$ the image of $(\mu;a_1,\dots,a_n)\in \P\circ A$ by the structural evaluation morphism.

  If $\P$ is a reduced operad, the endofunctor $\Gamma(\P)$ is a monad, and we define divided power $\P$-algebra, or $\Gamma(\P)$-algebras, to be the algebras over the monad $\Gamma(\P)$. The multiplication of this monad is given, for all vector spaces $V$, by the composite:
  \[
    \diag@=1.2cm{\P\tcirc\left(\P\tcirc V\right)\ar[r]& (\P\tcirc\P)\tcirc V\ar[r]^-{Tr^{-1}_{\P,\P}\tcirc V}& (\P\circ\P)\tcirc V\ar[r]^-{\mu_\P\circ V}& \P\tcirc V,}
  \]
  where the first arrow is given by the associator.

  In the case where $\P$ is a reduced operad, the natural transformation $Tr_\P:S(\P)\to \Gamma(\P)$ is a morphism of monads. When the characteristic of $\F$ is 0, $Tr_\P$ is an isomorphism.
  \end{defi}
  In Section \ref{secchardiv}, we will recall the characterisation of \cite{SI} for divided power algebras over an operad $\P$ in terms of monomial operations and relations (Theorem \ref{theoinv}) and introduce our notation for elements of the vector space $\Gamma(\M,V)$ (see Notation \ref{notetcirc}). This notation uses the product $\diamond$ for partition of integers which is defined in \ref{diamond}. We will give the formulae for the trace maps and for the associator of the product $\tcirc$ under this notation (See Lemmas \ref{lemtr} and \ref{lemass}).

  The trace map and associator are the principal ingredients to understanding interactions between the functor $\Gamma$ and the composition product of operads. They will play a key role in Section \ref{remdistrgamma} to obtain our main result on divided power algebras over a product of operads with distributive law.
  \begin{defi}[The operads $\As$, $\Com$, $\Lie$, $\D$ and $\Lev_p$]\label{exop}\item
    \begin{enumerate}[label=\alph*)]
      \item The symmetric sequence $\As=(\As(n))_{n>0}$, where $\As(n)$ is the regular representation, is endowed with an operad structure such that $\As$-algebras are exactly (non-unital) associative algebras (see \cite{LV}).
      \item  The symmetric sequence $\Com=(\Com(n))_{n>0}$, where $\Com(n)$ is the trivial representation, is endowed with an operad structure such that $\Com$-algebras are exactly (non-unital) associative and commutative algebras (see \cite{LV}). We denote by $X_n$ a generator of $\Com(n)$, so that $X_1=1_{\Com}$, and $X_n\circ_i X_m=X_{n+m-1}$.
      \item There is a symmetric sequence $\Lie=(\Lie(n))_{n>0}$, endowed with an operad structure, such that $\Lie$-algebras are exactly Lie algebras (see \cite{LV}). As an operad, $\Lie$ is generated by an element $[-;-]\in\Lie(2)$ such that $(\tau)\cdot [-;-]=-[-;-]$ and 
    \[
      [-;-]\circ_1[-;-]+(1\,2\,3)\cdot[-;-]\circ_1[-;-]+(1\,3\,2)\cdot[-;-]\circ_1[-;-]=0,
    \]
    where $\tau\in\∑_2$ is the transposition and we use the standard notation of the cycles $(1\,2\,3)$ and $(1\,3\,2)\in\∑_3$.
      \item The polynomial algebra $\D=\F[d]$ on one generator $d$, with its associative commutative and unital multiplication, can be seen as an operad concentrated in arity 1. $\D$-algebras are exactly vector spaces $A$ endowed with an endomorphism $d:A\to A$ (see \cite{SI3}).
      \item  For all integers $p$, there is an operad $\Lev_p$ such that a $\Lev_p$-algebra is a vector space $A$ endowed with a $p$-ary operation $\star$ such that, for all $a_1,\dots,a_p\in A$ and $\sigma\in\∑_p$,
    \[
      \star(a_1,\dots,a_p)=\star(a_{\sigma(1)},\dots,a_{\sigma(p)}),
    \]
    and such that for all $(a_{i,j})_{i\in[p],j\in[p]}\in A^{p^2}$,
    \begin{multline*}
      \star(\star(a_{1,1},a_{1,2},\dots,a_{1,p}),\dots,\star(a_{p,1},a_{p,2},\dots,a_{p,p}))\\=\star(\star(a_{1,1},a_{2,1},\dots,a_{p,1}),\dots,\star(a_{1,p},a_{2,p},\dots,a_{p,p})).
    \end{multline*}
    $\Lev_p$-algebras are called $p$-level algebras (see \cite{SI3}).
    \end{enumerate}
  \end{defi}
\section{Partitions and compositions of integers}\label{secpartcomp}
In this paper, a partition of a positive integer $n$ is an ordered partition of a finite set with $n$ elements. This convention differs from the classical terminology (see, for example, \cite{LEu}), and this choice of convention is justified by our use of the term partition. A composition of an integer $n$ is a finite ordered list of (positive) integers which sum up to $n$. Our definition differs from the classical terminology (see \cite{PAMM}) only by the fact that we authorise copies of 0 in our compositions.

This section, which is analogous to Section 2.1 of \cite{SI}, is devoted to the definition of several operations on partitions and compositions of integers. The three operations denoted by $\circ_i$, $\rhd$, and $\diamond$ will be used in Section \ref{secchardiv} to recall the characterisation of divided power algebras over an operad obtained in \cite{SI}. The operation $\diamond$ will also be used frequently in the sequel to describe the associator of the operadic product $\tcirc$ defined in Section \ref{modop}.
	\begin{Note*}\nopagebreak\item
  \begin{itemize}
  \item Denote by $[n]$ the set $\{1,\dots,n\}$.

  \item For a set $E$, $\∑_E$ is the group of bijections $E\to E$.

  \item For $R=(R_i)_{1\le i\le p}$ a partition of $[n]$, define $\∑_R:=\prod_{i=1}^p\∑_{R_i}\subset\∑_n$. 

  \item For $X$ an $\∑_n$-set, for $x\in X$ and for $R$ a partition of $[n]$, denote by $[x]_R$, and $Stab_R(x)$, the class and the stabiliser of $x$ under the action of the subgroup $\∑_R$ of $\∑_n$. 

  \item For $G$ a group, $X$ a $G$-set and $x\in X$, denote by $\Omega_G(x)=\{g\cdot x\ |\ g\in G\}$ the orbit of $x$ under the action of $G$.

  \item For $E\subseteq[n]$ and $\rho\in\∑_n$, denote by $\rho(E)$ the set $\{\rho(x):x\in E\}$. 

  \item For $E\subseteq\N$ and $n\in\N$, denote $E+n:=\{x+n:x\in E\}$.
  \end{itemize}
\end{Note*}
An ordered partition of the set $[n]=\{1,\dots,n\}$ is a finite sequence of (possibly empty) sets $J=(J_i)_{0\le i\le s}$ such that $\bigcup_{i=0}^sJ_i=[n]$, with $J_{i}\cap J_{i'}=\emptyset$ for all $i\neq i'$, where $s$ is any non-negative integer.
\begin{defi}\label{defiPi}
  There is an operad $\Pi$, called the operad of partitions, such that, for all $n\in\N$, $\Pi(n)$ is the vector space spanned by the ordered partitions $J=(J_i)_{0< i\le s}$ of the set $[n]$, where the length of the partition, $s$, can be arbitrarily large. The unit of the operad sends $1\in\F$ to the partition $(\{1\})\in\Pi(1)$ of $[1]$, and the compositions are induced by:
  \[(J\circ_l K)_i=\begin{cases}
      \lambda_{l,k}(J_i),&\mbox{if }l\in J_{i'}\mbox{ with }i'>i,\\
      \lambda_{l,k}(J_i)\setminus\{l\}\cup (K_{0}+l)&\mbox{if }l\in J_{i},\\
      \lambda_{l,k}(J_i)\cup (K_{i-i'}+l),&\mbox{if }l\in J_{i'}\mbox{ with }i'<i,
    \end{cases}\]
  where $K$ is a partition of $[k]$, $J$ is a partition of $[j]$, $l\in[j]$, and $\lambda_{l,k}\colon[j]\to[j-1+k]$ is the identity on $[l]$ and sends $l'>l$ to~$l'+k-1$.

  The symmetric group $\∑_n$ acts on $\Pi(n)$ by $\rho(R)=(\rho(R_i))_{i\in[p]}$.

  We also denote by $\bar\Pi$ the sub-operad of $\Pi$ such that $\bar\Pi(0)=0$ and $\bar\Pi(n)=\Pi(n)$ for all $n>0$.
\end{defi}
\begin{note}\item\nopagebreak
  \begin{itemize}\label{COMP}
    \item A $p$-tuple of non-negative integers $\r=(r_1,\dots,r_p)$ such that $r_1+\dots+r_p=n$ is called a composition of the integer $n$. Denote by $\Comp_p(n)$ the set of compositions of $n$ into $p$ parts and $\Pi_p(n)$ the set of partitions of the set $[n]$ into $p$ parts. There is an injection $\iota:\Comp_p(n)\inj \Pi_p(n)$, mapping $\r$ to:
    \[
    	r_1+\dots+r_{i-1}+[r_i]=\{r_1+\dots+r_{i-1}+1,\dots, r_1+\dots+r_i\}.
    \]
    The composition $\r$ will be identified with its image $\iota(\r)$. the map $\iota$ admits a left inverse $pr:\Pi_p(n)\surj\Comp_p(n)$, which sends $R=(R_i)_{i\in[p]}$ to $(\val{R_i})_{i\in[p]}$.
    \item $\Comp_p(n)$ is endowed with the following $\∑_p$-action: if $\sigma\in\∑_p$ and $\r\in\Comp_p(n)$, $\r^\sigma:=(r_{\sigma^{-1}(1)},\dots,r_{\sigma^{-1}(p)})$. Note that $\iota$ is not compatible with this action.
    \item If $\r\in \Comp_p(n)$ and $\_q\in \Comp_s(r_i)$, then $\r\circ_i\_q=\_t\in \Comp_{p-1+s}(n)$, with
    \[
    	\_t=(r_1,\dots,r_{i-1},q_1,\dots,q_s,r_{i+1},\dots,r_p).
    \]
    \end{itemize}
\end{note}
\begin{rema}\label{partfonc}
  Equivalently, the set $\Pi_p(n)$ can be defined as the set of maps $f:[n]\to[p]$. The identification is done as follows: a partition $R=(R_i)_{i\in[p]}$ of $[n]$ is identified with the function $\p_R:[n]\to[p]$ mapping $x\in[n]$ to the unique $i\in[p]$ such that $x\in R_i$. Conversely, each $f:[n]\to[p]$ is associated with the partition $(f^{-1}(i))_{i\in [p]}$ of $[n]$. For $\sigma\in\∑_p$ and $f:[n]\to[p]$, one has $\sigma\cdot f:x\mapsto \sigma(f(x))$.

  With this setting, for $f\in\Pi_p(n)$ and $g\in\Pi_s\left(|f^{-1}(i)|\right)$, there exists a unique increasing bijection $b:f^{-1}(i)\to \left[|f^{-1}(i)|\right]$, so that
    \begin{eqnarray*}
      f\circ_ig:&[n]&\to[p+s-1]\\
      &x&\mapsto\left\{
      \begin{array}{lcc}
      f(x),&\mbox{if}&f(x)<i,\\
      g\circ b(x)+i-1,&\mbox{if}&f(x)=i,\\
      f(x)+s-1,&\mbox{if}&f(x)\ge i.\\
      \end{array}\right.
    \end{eqnarray*}
\end{rema}
\begin{defis}\textbf{Operations on partitions}\label{wedge}\label{diamond}\\
  In the next paragraphs we will define four operations on partitions: 1) a product $\rhd$, 2) a product $\otimes$, 3) a family of unitary operations $(\gamma_k)_{k\in\N}$, and 4) a composition $\diamond$.
\begin{enumerate}[label=\arabic*)]
  \item Let $Q\in\Pi_s(p)$ and $R\in \Pi_p(n)$. According to the previous remark, $Q$ corresponds to $\p_Q:[p]\to[s]$ and $R$ to $\p_R:[n]\to[p]$. Define the partition
  \[
  	Q\rhd R\in\Pi_s(n)
  \]
  associated with the function $\p_Q\circ\p_R$, that is:
  \[
  	(Q\rhd R)_{i\in[s]}=\left(\coprod_{j\in Q_i}R_j\right)_{i\in[s]}.
  \]
Note that, if $\_q\in \Comp_s(p)$ and $\r\in \Comp_p(n)$, then $\_q\rhd\r=\_t\in \Comp_s(n)$, where $\_t=\big(\sum_{i\in \_q_1}r_i,\dots,\sum_{i\in\_q_s}r_i\big)$.

  \item The operation:
  \begin{eqnarray*}
    \otimes:&\Pi_p(n)\times\Pi_s(m)&\to\Pi_{p+s}(n+m)\\
    &((R_i)_{i\in[p]},(Q_j)_{j\in[s]})&\mapsto R\otimes Q:=({(n,m)}\circ_2Q)\circ_1R.
  \end{eqnarray*}
satisfies $R\otimes Q=(R_1,\dots,R_p,Q_1+n,\dots,Q_s+n)$ and is associative.
  \item For $k$ a non-negative integer, and $R\in\Pi_p(n)$, $R^{\otimes k}$ denotes the partition:
\[
	R^{\otimes k}:=\underset{k}{\underbrace{R\otimes\dots\otimes R}}\in\Pi_{kp}(kn).
\]
  For all positive integers $k$ and $R\in\Pi_p(n)$, set:
  \[
  	\gamma_k(R):=\left(\{i+(j-1)p\}_{j\in[k]}\right)_{i\in[p]}\rhd R^{\otimes k}=\left(\coprod_{j=0}^{k-1}R_i+jn\right)_{i\in[p]}\in\Pi_p(kn).
  \]
  
\begin{ex*}
  For $n=3$ and $p=2$, $R=(\{1,3\},\{2\})$, $k=3$,
  \[
  	R^{\otimes 3}=(\{1,3\},\{2\},\{4,6\},\{5\},\{7,9\},\{8\})\in\Pi_6(9),
  \]
  and:
  \[
  	\gamma_3(R)=(\{1,3,4,6,7,9\},\{2,5,8\})\in\Pi_2(9).
  \]
\end{ex*}

  \item Let $\r\in\Comp_p(n)$ and, for all $i\in[p]$, let $\q_i\in \Comp_{s_i}(m_i)$. For all $Q_1\in\Pi_{s_1}(m_1)$, $Q_2\in\Pi_{s_2}(m_2)$,\dots,$Q_p\in\Pi_{s_p}(m_p)$, set:
  \[
  	\r\diamond(Q_1,\dots,Q_p)=\gamma_{r_1}(Q_1)\otimes\dots\otimes\gamma_{r_p}(Q_p).
  \]
\begin{ex*}
  For $n=5$, $p=2$, $\r=(3,2)$, $m_1=m_2=3$, $s_1=s_2=2$, $\q_1=(2,1)$ and $\q_2=(1,2)$,
  \begin{align*}
    \r\diamond(\q_1,\q_2)&=\gamma_3(\q_1)\otimes\gamma_2(\q_2)\\
    &=(\{1,2,4,5,7,8\},\{3,6,9\})\otimes (\{1,4\},\{2,3,5,6\})\\
    &=(\{1,2,4,5,7,8\},\{3,6,9\},\{10,13\},\{11,12,14,15\})\in\Pi_4(15).
  \end{align*}
\end{ex*}
\end{enumerate}
\end{defis}
	\section{Recollections about distributive laws}\label{secDL}
	 For two monads $T,S$ on the same base category, there is no natural way \textit{a priori} to define a monad structure on the composition product $T\circ S$. Indeed, such a monad structure would include a unit - for which a natural choice is the product of the two units of $T$ and $S$, and a composition map 
	 \[
	 	\mu:\left(T\circ S\right)\circ \left(T\circ S\right)\to T\circ S.
	 \]
	 There is no natural choice for such a map $\mu$. A distribution law is a natural transformation $S\circ T\to T\circ S$ which can be used to naturally define such a composition $\mu$.

	 In this section, we introduce this notion of distributive law on a product of monads due to Beck \cite{JBe}, and the analogue notion of operadic distributive law on a product of operads \cite{FM,MMa,SML}. We introduce the classical example of the operad $\Pois$, which is obtained from the product of operads $\Com\circ\Lie$ \textit{via} a distributive law. We also state a result due to Beck \cite{JBe} which shows that, given two monads $T$ and $S$, any monadic structure on $T\circ S$ compatible with the respective units of $T$ and $S$, can be obtained from a distributive law (see Proposition \ref{propbeck}). This result, and the sketch of proof given in Remark \ref{propbeckrema}, will play the main role in identifying divided power structures on product of operads with distributive laws (see Section \ref{secdistrlawdiv}).
	\begin{defi}[\cite{SML}]\label{defDL}
		Let $S,T:\C\to \C$ be two monads in $\C$. A distributive law is a natural transformation
		\[
			\Lambda:T\circ S\to S\circ T
		\]
		such that the following diagrams commute:

		\begin{minipage}[t]{.53\linewidth}
			\begin{equation}\tag{DL1}\label{DL1}
			\diag@=.8cm{
				T\circ S\circ S\ar[r]^-{\Lambda\circ S}\ar[d]_-{ T\circ\mu_{S}}&S\circ T\circ S\ar[r]^-{S\circ\Lambda}&S\circ S\circ  T\ar[d]^-{\mu_S\circ T}\\
				T\circ S\ar[rr]_-{\Lambda}&&S\circ T,
			}
		\end{equation}
		\end{minipage}
		\begin{minipage}[t]{.42\linewidth}
			\begin{equation}\tag{DL2}\label{DL2}
			\diag@=.8cm{
			& T\ar[dl]_-{ T\circ\eta_{S}}\ar[rd]^-{\eta_{S}\circ T}\\
			 T\circ S\ar[rr]^-{\Lambda}&&S\circ T,
			}
		\end{equation}
		\end{minipage}

		\begin{minipage}[t]{.53\linewidth}
			\begin{equation}\tag{DL3}\label{DL3}
			\diag@=.8cm{
				T\circ  T\circ S\ar[r]^-{ T\circ\Lambda}\ar[d]_-{\mu_{T}\circ S}& T\circ S\circ T\ar[r]^-{\Lambda\circ T}&S\circ  T\circ T\ar[d]^-{S\circ\mu_{ T}}\\
				 T\circ S\ar[rr]_-{\Lambda}&&S\circ T,
			}
		\end{equation}
		\end{minipage}
		\begin{minipage}[t]{.42\linewidth}
			\begin{equation}\tag{DL4}\label{DL4}
			\diag@=.8cm{
			&S\ar[dl]_-{\eta_{T}\circ S}\ar[rd]^-{S\circ\eta_{T}}\\
			 T\circ S\ar[rr]^-{\Lambda}&&S\circ T.
			}
		\end{equation}
		\end{minipage}
	\end{defi}
	\begin{prop}[\cite{JBe},\cite{SML}]\label{propdistrmonad}
		If $\Lambda:T\circ S\to S\circ T$ is a distributive law, then $S\circ T$ is endowed with a monad structure in $\C$ such that $\eta_{S\circ T}=\eta_S\circ\eta_{T}$ and $\mu_{S\circ T}=(\mu_S\circ\mu_T)(S\circ\Lambda\circ T)$.
	\end{prop}
  The following definition gives conditions for a map of operads $\Lambda:\Q\circ\P\to\P\circ\Q$ to induce a distributive law on their associated Schur functors.
	\begin{defi}[\cite{LV}]\label{defODL}
		Let $\P,\Q$ be two operads. An operadic distributive law is a morphism of symmetric sequences:
		\[
			\Lambda:\Q\circ\P\to \P\circ\Q
		\]
		such that the following diagrams commute:

		\begin{minipage}[t]{.6\linewidth}
			\begin{equation}\tag{ODL1}\label{ODL1}
			\diag@=.8cm{
				\Q\circ\P\circ\P\ar[r]^-{\Lambda\circ\P}\ar[d]_-{\id_\Q\circ\mu_\P}&\P\circ\Q\circ\P\ar[r]^-{\P\circ\Lambda}&\P\circ\P\circ\Q\ar[d]^-{\mu_\P\circ\id_\Q}\\
				\Q\circ\P\ar[rr]_-{\Lambda}&&\P\circ\Q,
			}
		\end{equation}
		\end{minipage}
		\begin{minipage}[t]{.35\linewidth}
			\begin{equation}\tag{ODL2}\label{ODL2}
			\diag@=.45cm{
			& \Q\ar[dl]_-{ \Q\circ\eta_{\P}}\ar[rd]^-{\eta_{\P}\circ \Q}\\
			 \Q\circ \P\ar[rr]^-{\Lambda}&&\P\circ \Q,
			}
		\end{equation}
		\end{minipage}
		
		\begin{minipage}[t]{.6\linewidth}
			\begin{equation}\tag{ODL3}\label{ODL3}
			\diag@=.8cm{
				\Q\circ  \Q\circ \P\ar[r]^-{ \Q\circ\Lambda}\ar[d]_-{\mu_{\Q}\circ \P}& \Q\circ \P\circ \Q\ar[r]^-{\Lambda\circ \Q}&\P\circ  \Q\circ \Q\ar[d]^-{\P\circ\mu_{ \Q}}\\
				 \Q\circ \P\ar[rr]_-{\Lambda}&&\P\circ \Q,
			}
		\end{equation}
		\end{minipage}
		\begin{minipage}[t]{.35\linewidth}
			\begin{equation}\tag{ODL4}\label{ODL4}
			\diag@=.41cm{
			&\P\ar[dl]_-{\eta_{\Q}\circ \P}\ar[rd]^-{\P\circ\eta_{\Q}}\\
			 \Q\circ \P\ar[rr]^-{\Lambda}&&\P\circ \Q.
			}
		\end{equation}
		\end{minipage}
	\end{defi}
	\begin{prop}[\cite{LV}]
		If $\Lambda:\Q\circ\P\to\P\circ\Q$ is an operadic distributive law, then $\P\circ\Q$ is endowed with the structure of an operad such that $\eta_{\P\circ\Q}=\eta_\P\circ\eta_\Q$ and $\mu_{\P\circ\Q}=(\mu_\P\circ\mu_\Q)(\P\circ\Lambda\circ\Q)$, and in this case, $S(\Lambda):S(\Q\circ\P)=S(\Q)\circ S(\P)\to S(\P)\circ S(\Q)$ is a distributive law, and the resulting monad structure on $S(\P)\circ S(\Q)$ corresponds to the monad structure of $S(\P\circ \Q)$ induced by the operad structure of $\P\circ\Q$.
	\end{prop}
  \begin{ex}[\cite{BF}]\label{exPois}
    There is a distributive law $\rho:\Lie\circ\Com\to\Com\circ\Lie$ induced by:
    \[
     \rho([-;-]; X_2,1_{\Com})= (2\, 3)\cdot(X_2;1_{\Lie},[-;-])+(X_2;[-,-],1_{\Lie}).
    \]
    The resulting operad is denoted $\Pois$, and $\Pois$-algebras are exactly the classical Poisson algebras.
  \end{ex}
  \begin{rema}
    Bremner--Markl \cite{BrMa} and Bremner--Dotsenko \cite{BD} have classified the distributive laws between the operads $\As$, $\Com$ and $\Lie$.
  \end{rema}
	\begin{prop}[\cite{JBe}]\label{propbeck}
		Let $S,T:\C\to\C$ be two monads in $\C$. There is a bijective correspondence between distributive laws $T\circ S\to S\circ T$ and maps $m:S\circ T\circ S\circ T\to S\circ T$ such that $S\circ T$, with $m$ and the unit $\eta_S\circ\eta_T$, is a monad.
	\end{prop}
	\begin{rema}[\cite{JBe}]\label{propbeckrema}
		One direction of the correspondence of Proposition \ref{propbeck} is given by Proposition \ref{propdistrmonad}. For the other direction, if $(S\circ T,m,\eta_S\circ\eta_T)$ is a monad, then $m\circ(\eta_S\circ T\circ S\circ \eta_T):T\circ S\to S\circ T$ is the desired distributive law.
	\end{rema}
  We will use the following algebraic fact:
	\begin{fact}\label{transfmonoid}
		Let $L$ be an object of a monoidal category $(\mathcal C,\otimes)$, let $(M,\mu,\eta)$ be a monoid object in $\mathcal C$, and let $\phi:L\to M$ be an isomorphism. The object $L$ is endowed with a monoid structure with unit $\phi^{-1}\circ \eta$ and multiplication $\phi^{-1}\circ\mu\circ(\phi\otimes\phi)$.
	\end{fact}
	\section{Shifted and derivation algebras}\label{secshiftder}
	In Section \ref{modop}, we introduced an operad $\D$ whose algebras are a vector spaces $A$ endowed with a linear map $d:A\to A$. For any operad $\P$, a distributive law $\Lambda$ on the product $\P\circ\D$ induces a category of algebras whose objects are $\P$-algebras $A$ endowed with a linear map $d:A\to A$, such that $d$ satisfies some compatibility relations with the operations of $\P$, these relations being induced by $\Lambda$. 

	In this section, we define two different distributive laws $\lambda$ and $\rho$ on $\P\circ \D$. The first one encodes the relation for $d$ to be a $\P$-endomorphism. The second one encodes the relation for $d$ to be a $\P$-derivation, in the sense of \cite{MSS}. This operad $\P\circ\D$, with the distributive law $\rho$, corresponds to the operad of $\P$-algebras with derivation of \cite{JLLo}.
	\begin{defi}[\cite{SI3}]
		For all operads $\P$, there is an operad $\Shift_{\P}$ obtained from the distributive law $\lambda:\D\circ\P\to\P\circ\D$ sending the generic element $(d^j;\mu)$ to $(\mu;d^j,\dots,d^j)$.

		A $\Shift_\P$-algebra is a $\P$-algebra equipped with a $\P$-algebra endomorphism $d$.
	\end{defi}
  \begin{prop}[\cite{SI3}]\label{expartcom}
    There is an isomorphism between $\Shift_{\Com}$ and the sub-operad $\bar\Pi$ of the operad of partition $\Pi$ from Definition \ref{defiPi}.
  \end{prop}
  \begin{defi}[\cite{MSS}]
    Let $\P$ be an operad and $A$ be a $\P$-algebra. A $\P$-derivation of $A$ is a linear map $d:A\to A$ satisfying, for all $\mu\in\P(n)$, and $a_1,\dots,a_n\in A$,
    \[
      d\left(\mu(a_1,\dots,a_n)\right)=\sum_{i=1}^{n}\mu(a_1,\dots,a_{i-1},da_i,a_{i+1},\dots,a_n).
    \]
  \end{defi}
	\begin{prop}
		The morphism of symmetric sequences $\rho:\D\circ\P\to\P\circ\D$ sending $(d^j;\mu)$ to $\sum_{\q\in\Comp_m(j)}\binom{j}{q_1,\dots,q_m}(\mu;d^{q_1},\dots,d^{q_m})$, where $\mu\in\P(m)$, is an operadic distributive law. The resulting operad will be denoted by $\Der_\P$.

		A $\Der_\P$-algebra is a $\P$-algebra equipped with a $\P$-derivation $d$.
	\end{prop}
	\begin{proof}
		We have to prove that the diagrams \ref{ODL1} to \ref{ODL4} of Definition \ref{defODL} commute, with $\Q=\D$, and $\Lambda=\rho$. Checking that diagrams \ref{ODL2} and \ref{ODL4} is a straightforward verification. Since $\D$ is generated by $d$, to prove that diagram \ref{ODL1} commutes, it suffices to prove that the diagram commutes on elements of the type $(d;\mu;\nu_1,\dots,\nu_k)\in\D\circ\P\circ\P$, with $\nu_i\in\P(n_i)$ for all $i\in[m]$. This is again a straightforward verification.
		
		
		We will prove that diagram \ref{ODL3} commutes on generic elements $(d^j;d^k;\mu)\in\D\circ\D\circ\P$. Following the diagram \ref{ODL3} counter-clockwise, one gets the elements $(\mu_D\circ\P)(d^j;d^k;\mu)=(d^{j+k};\mu)$, and
		\[
			\Lambda(d^{j+k};\mu)=\sum_{\q\in\Comp_m(j+k)}\binom{j+k}{q_1,\dots,q_m}(\mu,d^{q_1},\dots,d^{q_m}).
		\]
		 Following the diagram \ref{ODL3} clockwise, one gets the element:
		\[
			(\D\circ\Lambda)(d^j;d^k;\mu)=\sum_{\q''\in\Comp_m(k)}\binom{k}{q''_1,\dots,q''_m}(d^j;\mu;d^{q''_1},\dots,d^{q''_m}).
		\]
      Then, applying $(\Lambda\circ\D)$ to this element, one gets:
		\[
			\sum_{\q'\in\Comp_m(j),\q''\in\Comp_m(k)}\binom{j}{q_1',\dots,q_m'}\binom{k}{q''_1,\dots,q''_m}\left(\mu;d^{q'_1},\dots,d^{q'_m};d^{q''_1},\dots,d^{q''_m}\right),
		\]
		and finally, $(\P\circ\mu_{\D})\circ(\Lambda\circ\D)\circ(\D\circ\Lambda)(d^j;d^k;\mu)$ is equal to:
		\[
			\sum_{\q'\in\Comp_m(j),\q''\in\Comp_m(k)}\binom{j}{q_1',\dots,q_m'}\binom{k}{q''_1,\dots,q''_m}\left(\mu;d^{q''_1+q'_1},\dots,d^{q''_m+q'_m}\right),
		\]
		which is equal to:
		\[
			\sum_{\q\in\Comp_m(j+k)}\left(\sum_{(q',q")\in E_\q}\binom{j}{q'_1,\dots,q'_m}\binom{k}{q''_1,\dots,q''_m}\right)\left(\mu;d^{q_1},\dots,d^{q_m}\right),
		\]
		where $E_\q$ is the subset of $\Comp_m(j)\times\Comp_m(k)$ formed by pairs $(\q',\q'')$ satisfying, for all $i\in[m]$, $q'_i+q''_i=q_i$. To prove the equality:
 		\[
 		 	\left(\sum_{(q',q'')\in E_\q}\binom{j}{q'_1,\dots,q'_m}\binom{k}{q''_1,\dots,q''_m}\right)=\binom{j+k}{q_1,\dots,q_m},
 		 \]
 		which is a variant of the classical Vandermonde identity on binomial coefficients, note that in the polynomial algebra $\F[x_1,\dots,x_{m-1}]$, one has, on one hand,
 		\[
 			(1+x_1+\dots+x_{m-1})^{j+k}=\sum_{\q\in\Comp_m(j+k)}\binom{j+k}{q_1,\dots,q_m}x_1^{q_1}\dots x_{m-1}^{q_m},
 		\]
 		and on the other hand,
 		\begin{multline*}
 			(1+x_1+\dots+x_{m-1})^{j+k}=(1+x_1+\dots+x_{m-1})^{j}(1+x_1+\dots+x_{m-1})^{k},\\
 			=\left(\sum_{\q'\in\Comp_m(j)}\binom{j}{q'_1,\dots,q'_m}x_1^{q_1}\dots x_{m-1}^{q_{m-1}}\right)\left(\sum_{\q''\in\Comp_m(k)}\binom{k}{q''_1,\dots,q''_m}x_1^{q''_1}\dots x_{m-1}^{q_{m-1}}\right),\\
 			=\sum_{\q\in\Comp_m(j+k)}\left(\sum_{(q',q'')\in E_\q}\binom{j}{q'_1,\dots,q'_m}\binom{k}{q''_1,\dots,q''_m}\right)x_1^{q_1}\dots x_{m-1}^{q_{m-1}}.
 		\end{multline*}
	\end{proof}
  \begin{rema}
    In the preceding proof, when showing that diagram \ref{ODL3} commutes, we haven't used the fact that $\D$ was generated by $d$. Checking that this diagram commutes only on elements of the type $(d;d;\mu)$ would prove that our definition of $\rho$ on elements of the type $(d;\mu)$ induces a distributive law, but not that this distributive law coincides with the expression of $\rho$ that we have given for more general elements.
  \end{rema}

\begin{ex}
      Algebras over $\Der_{\Com}$ (resp. $\Der_{\Lie}$) are commutative associative algebras (resp. Lie algebras) equipped with a derivation in the classical sense.
\end{ex}

 \section{Characterisation of divided power structures}\label{secchardiv}
 In this section, we recall the main result of \cite{SI}, which states that there is an isomorphism of categories between the category of divided power algebras over an operad $\P$ and a category of vector spaces equipped with a family of monomial operations $(\beta_{x,\r})_{x,\r}$, indexed by elements of $\P$ and compositions of integers, satisfying a list of eight relations denoted by ($\beta$1) to ($\beta$7a) and ($\beta$7b). This result allows one to work on elements of a divided power algebra $A$ over an operad $\P$ without reference to the structural map $\Gamma(\P,A)\to A$. The operations $\beta_{x,\r}$ are to be considered as divided power operations, in the sense that they act on elements of the algebra as divided monomials (see Remark \ref{remaPtr}), following the classical example of divided $\Com$-algebras (see Example \ref{exdivop} and \cite{NR,BF}).

 In \ref{notetcirc} we introduce a new notation for elements of a product of symmetric sequences of the type $\M\tcirc\ı$, inspired by the notation used in Theorem \ref{theoinv}. Lemmas \ref{lemtr} and \ref{lemass} give the formulae for the trace map and the associator using this notation. These two technical results are the key to understanding Section \ref{secdistrlawdiv} and the compatibility between divided power algebras and distributive laws (see Proposition \ref{propdistrgamma} and Lemma \ref{lemmtLambda}).

 In Example \ref{exdivop}, we compare the result of Theorem \ref{theoinv} to the structure, characterised by Fresse \cite{BF}, of divided power algebras over the operads $\As$, $\Com$ and $\Lie$. These characterisations will be used in the following sections to identify divided power structures over products of these operads with a distributive law.
  \begin{theo}[\cite{SI}]\label{theoinv}
    A divided power algebra over a reduced operad $\P$ is a vector space $A$ endowed with a family of operations $\beta_{x,\r}:A^{\times p}\to A$, given for all $\r\in\Comp_p(n)$ and $x\in\M(n)^{\∑_\r}$, and which satisfy the relations:
  \begin{enumerate}[label=($\beta$\arabic*),itemsep=5pt]
    \item\label{relperm}$\beta_{x,\r}((a_i)_i)=\beta_{\rho^*\cdot x,\r^\rho}((a_{\rho^{-1}(i)})_{i})$ for all $\rho\in\∑_p$, where $\rho^*$ denotes the block permutation with blocks of size $(r_i)$ associated to $\rho$.

    \item\label{rel0} $\beta_{x,(0,r_1,r_2,\dots,r_p)}(a_0,a_1,\dots,a_p)=\beta_{x,(r_1,r_2,\dots,r_p)}(a_1,\dots,a_p)$.

    \item\label{rellambda} $\beta_{x,\r}(\lambda a_1,a_2,\dots,a_p)=\lambda^{r_1}\beta_{x,\r}(a_1,\dots,a_p)\quad\forall \lambda\in \F$.

    \item\label{relrepet}
    If $\r\in \Comp_p(n)$ and $\q\in \Comp_s(p)$, then
    \[
      \beta_{x,\r}(\underbrace{a_1,\dots,a_1}_{q_1},\underbrace{a_2,\dots,a_2}_{q_2},\dots,\underbrace{a_s,\dots,a_s}_{q_s})=\beta_{\big(\sum_{\sigma\in \∑_{\q\rhd\r}/\∑_\r}\sigma\cdot x\big),\ \q\rhd\r}(a_1,a_2,\dots,a_s).
    \]

    \item\label{relsomme} $\beta_{x,\r}(a_0+a_1,\dots,a_p)=\sum_{l+m=r_1}\beta_{x,\r\circ_1(l,m)}(a_0,a_1,\dots,a_p)$.

    \item\label{rellin} $\beta_{\lambda x+y,\r}=\lambda\beta_{x,\r}+\beta_{y,\r}$ , for all $x,y\in\M(n)^{\∑_\r}$,
  \end{enumerate}
   \begin{enumerate}[label=($\beta$7\alph*)]
    \item \label{relunit}$\beta_{1_\P,(1)}(a)=a \quad \forall a\in A$.
    \item \label{relcomp} Let $\r\in \Comp_p(n)$, $x\in \P(n)^{\∑_\r}$ and for all $i\in[p]$, let $\q_i\in \Comp_{k_i}(m_i)$, $x_i\in\P(m_i)^{\∑_{\q_i}}$ and $(b_{ij})_{1\le j\le k_i}\in A^{\times k_i}$. Denote by $b=(b_{ij})_{i\in[p],j\in[k_i]}$ and for all $i\in[p]$, $b_i=(b_{i,j})_{j\in[k_i]}$. Then:
    \[
      \beta_{x,\r}(\beta_{x_1,\q_1}(b_{1}),\dots,\beta_{x_p,\q_p}(b_{p}))=\beta_{\sum_{\tau}\tau\cdot\mu \big(x\otimes \big(\bigotimes_{i=1}^px_i^{\otimes r_i}\big)\big),\r\diamond(\q_i)_{i\in[p]}}(b),
    \]
    where $\r\diamond(\q_i)_{i\in[p]}$ is defined in \ref{diamond}, where $\beta_{\cdot,\r\diamond(\q_i)_{i\in[p]}}$ is defined in Remark \ref{rempart} and where $\tau$ ranges over $\∑_{\r\diamond(q_i)_{i\in[p]}}/(\prod_{i=1}^{p}\∑_{r_i}\wr \∑_{\q_i})$ in the sum.
  \end{enumerate}

  morphisms of divided power algebras over $\P$ are linear maps compatible with the operations $\beta_{x,\r}$.
  \end{theo}
  \begin{rema}\label{rempart}  
For a divided power $\P$-algebra $A$, the operations $\beta_{x,\r}$ induce operations $\beta_{x,R}:A^{\times p}\to A$, for all partitions $R\in\Pi_p(n)$ and any $x$ stable under the action of $\∑_R$. Indeed, there exists $\tau\in\∑_n$ such that $\tau(R)=\r$. Generally, there are several choices for the permutation $\tau\in\∑_n$ such that $\tau(R)=\r$. One can define:
\[
  \beta_{x,R}(a_1,\dots,a_p)=\beta_{\tau\cdot x,\r}(a_1,\dots,a_p),
\]
which does not depend on the chosen $\tau$.
\end{rema}
\begin{note}\label{notetcirc}
    Recall that, for all symmetric sequences $\M,\mathcal N$ such that $\M$ is reduced, $\M\tcirc \mathcal N$ is generated, as a vector space, by the elements $\t$ such that:
    \[
      \t=\sum_{\sigma\in\∑_n/\∑_\r}\sigma \left(x\otimes \bigotimes_{i\in[p]}y_i^{\otimes r_i}\right),
    \]
    with $x\in\M(n)^{\∑_\r}$, and $y_1\dots y_p\in\mathcal N$.

    The element $\t$ will be denoted by $\beta_{x,\r}(y_1,\dots,y_p)$. In the case where $\M=\P$ is an operad and $\mathcal N=A$ is a $\Gamma(\P)$-algebra concentrated in arity $0$, the element $\beta_{x,\r}(a_1,\dots,a_p)\in A$ is the image under the structural map $\P\tcirc A\to A$ of the element $\sum_{\sigma\in\∑_n/\∑_\r}\sigma \left(x\otimes \bigotimes_{i\in[p]}a_i^{\otimes r_i}\right)$, which justifies this notation \textit{a posteriori}.
  \end{note}
  \begin{lemm}\label{lemtr}
    The trace map $\M\circ\mathcal N\to \M\tcirc\mathcal N$ sends $(x;y_1,\dots,y_n)$ to $\beta_{x,(1,\dots,1)}(y_1,\dots,y_n)$.
  \end{lemm}
  \begin{proof}
    This is readily checked.
  \end{proof}
  \begin{lemm}\label{lemass}
     The associator $\as_{\M,\ı,\O}:\M\tcirc(\ı\tcirc \O)\to (\M\tcirc \ı)\tcirc \O$ sends
    \[
      \beta_{x,\r}(\beta_{y_1,\q_1}(z_{1,1},\dots,z_{1,s_1}),\dots,\beta_{y_p,\q_p}(z_{p,1},\dots,z_{p,s_p}))
    \]
    to
    \[
      \beta_{\beta_{\sum_{\sigma\in\∑_{\r\diamond(\q_1,\dots,\q_p)}/\prod_{i\in[p]}\∑_{r_1}\wr\∑_{\q_i}}x,(1,\dots,1)}(y_1^{\times r_1},\dots,y_p^{\times r_p}),\r\diamond(\q_1,\dots,\q_p)}(z_{1,1},\dots,z_{1,s_1},\dots,z_{p,s_p}).
    \]
  \end{lemm}
  \begin{proof}
    This is deduced from \cite{SI}, appendix A, using the fact that $\prod_{i\in[p]}\∑_{r_1}\wr\∑_{\q_i}$ is a subgroup of $\∑_{\r\diamond(\q_1,\dots,\q_p)}$.
  \end{proof}
  \begin{rema}\label{remaPtr}
   Let $A$ be a divided $\P$-algebra. Then $A$ is endowed with a structure of $\P$-algebra induced by restriction along the trace map $Tr_\P:S(\P)\to \Gamma(\P)$. Following Lemma \ref{lemtr}, an element $x\in\P(n)$ acts on $A$ \textit{via} the operation $\beta_{x,(1,\dots,1)}$ ($n$ copies of $1$ in the partition).

   Using relation \ref{relrepet}, one can check that, for all $\r\in\Comp_p(n)$, and assuming $x$ is fixed by the action of $\∑_{\r}$,
   \[
     \beta_{x,(1,\dots,1)}(\underbrace{a_1,\dots,a_1}_{r_1},\dots,\underbrace{a_p,\dots,a_p}_{r_p})=\left(\prod_{i=1}^pr_i!\right)\beta_{x,\r}(a_1,\dots,a_p).
   \]
   In particular, if $\prod_{i=1}^pr_i!$ does not divide the characteristic of the base field, the relation
   \[
     \frac{\beta_{x,(1,\dots,1)}(\underbrace{a_1,\dots,a_1}_{r_1},\dots,\underbrace{a_p,\dots,a_p}_{r_p})}{\left(\prod_{i=1}^pr_i!\right)}=\beta_{x,\r}(a_1,\dots,a_p)
   \]
   justifies calling the $\beta_{x,\r}$ `divided power operations'. This terminology was first used by Cartan \cite{HC}, in French `\textit{opérations aux puissances divisées}', to study a type of algebra which was later identified as the $\Gamma(\Com)$-algebras (see Example \ref{exdivop}, \cite{NR,BF}).
  \end{rema}
\begin{ex}[Divided power $\As$, $\Com$, $\Lie$, $\D$-algebras]\label{exdivop}\item
  \begin{enumerate}[label=\alph*)]
      \item A $\Gamma(\As)$-algebra is an $\As$-algebra (see \cite{BF}). Denote by $\mu\in\As(2)$ the generator of the operad $\As$. The non-unital associative multiplication is given by the operation $\beta_{\mu,(1,1)}$.
      \item  A $\Gamma(\Com)$-algebra is a divided power algebra in the classical sense (see \cite{BF}, \cite{NR}), that is, a vector space $A$ endowed with a commutative and associative multiplication, and for all $n>0$, with a set map $\gamma_n:A\to A$ satisfying a list of relations. In \cite{SI}, we show that the multiplication is given by the operation $\beta_{X_2,(1,1)}$, and that the operation $\gamma_n$ is given by $\beta_{X_n,(n)}$.
      \item Following Soublin \cite{JPSou}, if the characteristic of the base field $\F$ is $p$, then a $\Gamma(\Com)$-algebra is a $\Com$-algebra $A$ endowed with a divided $p$-power $\gamma_p$, that is, a set map $\gamma_p:A\to A$ satisfying the four relations:
\begin{enumerate}[label=(C$p$\arabic*)]
    \item\label{rel2Com1} $\gamma_p(\lambda a)=\lambda^p\gamma_p(a)$
    \item\label{rel2Com2} $a^p=0$
    \item\label{rel2Com3} $\gamma_p(a+b)=\gamma_p(b)+\sum_{i=1}^{p-1}\frac{a^ib^{p-i}}{i!(p-i)!}+\gamma_p(a)$
    \item\label{rel2Com4} $\gamma_p(ab)=0$.
  \end{enumerate}
    \item Let $p$ be the characteristic of the base field $\F$. Following \cite{BF}, a $\Gamma(\Lie)$-algebra is the same thing as a $p$-restricted algebra, that is, a vector space endowed with a Lie bracket $[-;-]$ and a Frobenius operator $x\mapsto F(x)$, satisfying, for all $x,y\in A$, $\lambda\in\F$:
\begin{enumerate}[label=(L\arabic*),itemsep=5pt]
  \item $\ad(F(x))=F\left(\ad(x)\right)$,
  \item $F(\lambda x)=\lambda^pF(x)$,
  \item $F(x+y)=F(x)+F(y)+\sum_{i=1}^{p-1}\frac{s_i(x,y)}{i}$,
\end{enumerate}
where $\ad$ is the usual adjoint operator: $\ad(x)(y)=[x;y]$, and $s_i(x,y)$ is the coefficient of $\lambda^{i-1}$ in the formal expression $\ad(\lambda x+y)^{p-1}(x)$ (see \cite{NJa}). If we denote by $[-;-]\in\Lie(2)$ the operadic generator, the Lie bracket of a $\Gamma(\Lie)$-algebra is given by the operation $\beta_{[-;-],(1,1)}$, and, following the computations of \cite{BF}, one can check that the Frobenius operator $F$ is given by $\beta_{L_p,(p)}$, where $L_p\in\Lie(p)$ is the element:
\[
  \sum_{\sigma}\underbrace{[-;-]\circ_1[-;-]\circ_1\dots\circ_1[-;-]}_{p-1}\cdot \sigma,
\]
where the sum runs over the $\sigma\in\∑_p$ such that $\sigma(1)=1$.
  \item The trace map yields an isomorphism $S(\D)\to \Gamma(\D)$, so a $\Gamma(\D)$-algebra is a vector space $A$ endowed with an endomorphism $d:A\to A$.
    \end{enumerate}  

\end{ex}
	\section{Distributive laws on divided power functors}\label{secdistrlawdiv}\label{remdistrgamma}

	This section is devoted to showing that, in th case where $\P\circ\Q$ is a product of operads endowed with a distributive law, the resulting monad $\Gamma(\P\circ\Q)$ can itself be obtained from the composition of $\Gamma(\P)$ and $\Gamma(\Q)$ via a distributive law. This is the object of Proposition \ref{propdistrgamma}. This result in itself is guaranteed by a proposition due to Beck (see \ref{propbeck}), the proof of which allows us to compute explicitly the resulting distributive law. This distributive law will allow us to retrieve relations between the monomial operations $\beta_{x,\r}$ defined in Section \ref{secchardiv} coming from the operad $\P$ and from the operad $\Q$. Using the diagrams defining this distributive law, we produce a computational result (Lemma \ref{lemmtLambda}), which will help us characterise explicitly some divided power structure over products of operads in the following sections.

	Let $\P,\Q$ be two reduced operads endowed with a distributive law $\Lambda:\Q\circ\P\to\P\circ\Q$.
  
  The operad $\P\circ\Q$ is reduced, and, according to Definition \ref{defiop}, the multiplication of the monad $\Gamma(\P\circ\Q)$ is given, for all vector spaces $V$, by the composite:
  \begin{equation}\tag{D1}\label{diagdistrmult}
     \diag{(\P\circ\Q)\tcirc \left((\P\circ\Q)\tcirc V\right)\ar[r]^{\alpha_1}&\left((\P\circ\Q)\tcirc(\P\circ\Q)\right)\tcirc V\ar[r]^{\alpha_2}&\left((\P\circ\Q)\circ(\P\circ\Q)\right)\tcirc V\ar[d]^{\alpha_3}\\&&(\P\circ\Q)\tcirc V,}
  \end{equation}
  where the first arrow, is given by the associator: $\alpha_1=\as_{\P\circ\Q,\P\circ\Q,V}$, where the second arrow is the inverse of the trace map, $\alpha_2=Tr^{-1}_{\P\circ\Q,\P\circ\Q}(V)$ (see Lemma \ref{lemmTriso}), and the last arrow is itself the composite 
  \[
    \alpha_3=\mu_{\P\circ\Q}(V)=(\mu_{\P}\circ\mu_{\Q})(\id_{\P}\circ\Lambda\circ\id_\Q)(V).
  \]

	According to Lemma \ref{lemmTriso}, the trace map $Tr_{\P,\Q}$ is an isomorphism, so one can build an isomorphism $\Gamma(\P\circ\Q)\to \Gamma(\P)\circ\Gamma(\Q)$ given, for all vector spaces $V$, by:
	\[
		\diag@=1.2cm{(\P\circ\Q)\tcirc V\ar[r]^-{Tr_{\P,\Q}\tcirc V}&(\P\tcirc\Q)\tcirc V\ar[r]&\P\tcirc(\Q\tcirc V)},
	\]
	where the last arrow is given by the associator $\as_{\P,\Q,V}^{-1}$.

  We deduce from the Fact \ref{transfmonoid} that $\Gamma(\P)\circ \Gamma(\Q)$ inherits a monad structure from the monad structure of $\Gamma(\P\circ\Q)$. Since this isomorphism sends $\Gamma(\eta_{\P})\circ\Gamma(\eta_{\Q})$ to $\Gamma(\eta_{\P\circ\Q})$, Proposition \ref{propbeck} ensures that the monad structure on $\Gamma(\P)\circ\Gamma(\Q)$ is obtained from the monad structures of $\Gamma(\P)$ and $\Gamma(\Q)$ through a distributive law $\tLambda:\Gamma(\Q)\circ\Gamma(\P)\to \Gamma(\P)\circ\Gamma(\Q)$. More precisely, this distributive law is given for all vector spaces $V$ by:
	\begin{equation}\label{tildeD}\tag{D2}
    \diag@=.6cm{\Q\tcirc(\P \tcirc V)\ar[d]^{\delta_1}\\
    \P\tcirc \left(\Q\tcirc(\P\tcirc (\Q\tcirc V)\right)\ar[r]^{\delta_2}
    &\left((\P\tcirc\Q)\tcirc(\P\tcirc\Q)\right)\tcirc V\ar[r]^{\delta_3}&\left((\P\circ\Q)\tcirc(\P\circ\Q)\right)\tcirc V\ar[d]^{\delta_4}\\\P\tcirc(\Q\tcirc V)&(\P\tcirc\Q)\tcirc V\ar[l]^{\delta_6}&(\P\circ\Q)\tcirc V\ar[l]^{\delta_5},}
  \end{equation}
		
	where the first arrow is given by the unity morphisms of $\P$ and $\Q$: 
	\[
		\delta_1=\eta_\P\tcirc(\Q\tcirc(\P\tcirc(\eta_\Q\tcirc V))):\Q\tcirc(\P \tcirc V)\to \P\tcirc \left(\Q\tcirc(\P\tcirc (\Q\tcirc V)\right),
	\] 
	the second arrow $\delta_2$ is given by three consecutive use of the associator,
  \[
    \diag@R=1cm@C=3cm{\P\tcirc \left(\Q\tcirc(\P\tcirc (\Q\tcirc V)\right)\ar[r]^{\as_{\P,\Q,\P\tcirc (\Q\tcirc V)}}&(\P\tcirc\Q)\tcirc(\P\tcirc(\Q\tcirc V))\ar[d]^{\P\tcirc\Q\circ(\as_{\P,\Q,V})}\\
    ((\P\tcirc\Q)\tcirc(\P\tcirc\Q))\tcirc V&(\P\tcirc\Q)\tcirc((\P\tcirc\Q)\tcirc V)\ar[l]^{\as_{\P\tcirc\Q,\P\tcirc\Q,V}},}
  \]
  the third arrow is given by the trace maps:
	\[
		\delta_3=\left(Tr_{\P,\Q}^{-1}\tcirc Tr_{\P,\Q}^{-1}\right)\tcirc V:\left((\P\tcirc\Q)\tcirc(\P\tcirc\Q)\right)\tcirc V\to\left((\P\circ\Q)\tcirc(\P\circ\Q)\right)\tcirc V,
	\]
	the fourth arrow corresponds to the composition of the two last arrows of diagram \ref{diagdistrmult}), $\delta_4=\alpha_3\circ\alpha_2$, the fifth arrow is the trace map, $\delta_5=Tr_{\P\circ\Q}\tcirc V$, and the last arrow is the associator, $\delta_6=\as_{\P,Q,V}^{-1}$. We deduce the following:

  \begin{prop}\label{propdistrgamma}
    Let $\P,\Q$ be two reduced operads and $\Lambda:\Q\circ\P\to\P\circ\Q$ be an operadic distributive law. A $\Gamma(\P\circ\Q)$-algebra is a vector space endowed with both a $\Gamma(\P)$-algebra and a $\Gamma(\Q)$-algebra structures, with compatibility relation induced by the distributive law $\tLambda$.
  \end{prop}
  Note that this does not give us a general expression for the distributive law $\tLambda$. Consider a generic element of the form 
  \[
    \mathfrak s=\beta_{x,\r}(\beta_{y_1,\q_1}(a_{1,1},\dots,a_{1,s_1}),\dots,\beta_{y_p,\q_p}(a_{p,1},\dots,a_{p,s_p}))\in\Q\tcirc(\P\tcirc V),
  \]
  where $\r\in\Comp_p(n)$, $x\in\Q(n)^{\∑_{\r}}$, and for all $i\in[p]$, $\q_i\in\Comp_{s_i}(m_i)$, $y_i\in\P(m_i)^{\∑_{\q_i}}$, and $(a_{i,j})_{i,j}\in V^{s_1+\dots+s_p}$. Following Diagram \ref{tildeD}, starting with the element $\mathfrak s$, the first three maps can be given a general expression, independent from the operads $\P$ and $\Q$. This observation provides the following result, which will greatly simplify our computations in the sequel:
  \begin{lemm}\label{lemmtLambda}
    The element $\tLambda(\mathfrak s)$ is equal to $\as_{\P,\Q,V}\circ(Tr_{\P,\Q}\tcirc V)$ applied to the element:
    \begin{equation*}
      \beta_{\sum_{\sigma\in\∑_{\r\diamond(\q_1,\dots,\q_p)}/\prod_{i\in[p]}\∑_{r_i}\wr\∑_{\q_i}} \Lambda(x;y_1^{\times r_1},\dots,y_p^{\times r_p})\cdot \sigma,\r\diamond(\q_1,\dots,\q_p)}\left((a_{i,j})_{i\in[p],j\in[s_i]}\right),
    \end{equation*}
    for all $\mathfrak s\in\Q\tcirc(\P\tcirc V)$ as above, and where $\as_{\P,\Q,V}:(\P\tcirc \Q)\tcirc V\to \P\tcirc(\Q\tcirc V)$ is the associator.
  \end{lemm}
   
  Next Sections are devoted to applications of this result to particular operads with distributive laws.
	\section{Divided power algebras with shift/derivation operator}\label{secdivshiftder}
  In Section \ref{secshiftder}, we built, from an operad $\P$, two operads $\Shift_\P$ and $\Der_\P$, both obtained from the product $\P\circ\D$, respectively with the operadic distributive laws $\lambda$, and~$\rho$. In this section, we apply the results of Section \ref{secdistrlawdiv} to the two operads $\Shift_{\P}$ and $\Der_{\P}$. We obtain, in Theorems \ref{theodivshift} and \ref{theodivder}, a characterisation of divided power algebras over both of these operads. This leads us to define a notion of divided $\P$-derivation (see Definition \ref{defdivder}), which generalises the ``power rule'' on derivations over a divided power algebra, due to Keigher--Pritchard \cite{KP} (see Corollary \ref{expartdivder}).

  Following Proposition \ref{propdistrgamma}, assuming that $\P$ is reduced, a $\Gamma(\Shift_\P)$-algebra (respectively, a $\Gamma(\Der_\P)$-algebra) is a vector space endowed with both structures of $\Gamma(\P)$-algebra and of $\Gamma(\D)$-algebra, with compatibility relation between those structures given by two distributive laws $\tlambda$ and $\trho$ induced by $\lambda$ and $\rho$. Our first result is the following:
	\begin{theo}\label{theodivshift}
		A divided power $\Shift_\P$-algebra is a divided power $\P$-algebra endowed with an endomorphism of divided power $\P$-algebra.
	\end{theo}
  \begin{coro}\label{expartdiv}
    Recall that, as an operad, the operad of (non-empty) partitions $\bar\Pi$ is isomorphic to $\Shift_{\Com}$ (see Proposition \ref{expartcom}). As a result of Theorem \ref{theodivshift}, a divided power $\bar\Pi$-algebra is a divided power algebra in the classical sense of the term, endowed with an endomorphism of divided power algebra.
  \end{coro}
  \begin{proof}
      Following Proposition \ref{propdistrgamma}, a $\Shift_\P$-algebra is a vector space $A$ endowed with the structures of a $\Gamma(\P)$-algebra and of a $\Gamma(\D)$-algebra. We have seen (see Example \ref{exdivop}) that the $\Gamma(\D)$-algebra structure is the data of an endomorphism $d:A\to A$. To understand the compatibility relations between $d$ and the $\Gamma(\P)$-algebra structure, we have to compute the image of $d(\beta_{x,\r}(a_1,\dots,a_p))$ under the distributive law $\tlambda:\D\tcirc(\P\tcirc A)\to \P\tcirc(\D\tcirc A)$, for $x\in\P(n)^{\∑_\r}$. Applying Lemma \ref{lemmtLambda} to the element
      \[
        \mathfrak s=d(\beta_{x,\r}(a_1,\dots,a_p))=\beta_{d,(1)}(\beta_{x,\r}(a_1,\dots,a_p)),
      \]
      we obtain:
      \[
        \tlambda(\mathfrak s)=\as_{\P,\D,A}\circ (Tr_{\P,\D}\tcirc A)\left(\beta_{\sum_{\sigma\in\∑_{(1)\diamond \r}/\∑_{(1)}\wr\∑_{\r}} \lambda(d;x)\cdot\sigma,(1)\diamond(\r)}(a_1,\dots,a_p)\right).
      \]
      Note that $(1)\diamond \r=\r$, and $\∑_{(1)\diamond \r}=\∑_{(1)}\wr\∑_{\r}=\∑_{\r}$. So,
      \begin{align*}
        \tlambda(\mathfrak s)&=\as_{\P,\D,A}\circ(Tr_{\P,\D}\tcirc A)\left(\beta_{(x;d^{\times n}),\r}(a_1,\dots,a_p)\right).\\
        &=\as_{\P,\D,A}\left(\beta_{\beta_{x,(1,\dots,1)}(d,\dots,d),\r}(a_1,\dots,a_p)\right),\\
        &=\beta_{x,\r}(da_1,\dots,da_p),
      \end{align*}
      which proves that $d$ is an endomorphism of $\Gamma(\P)$-algebra.
  \end{proof}
  Understanding the expression of the morphism $\trho$ yields a notion of divided derivation:
  \begin{defi}\label{defdivder}
    Let $A$ be a divided power $\P$-algebra. A divided $\P$-derivation is a linear map $d:A\to A$ satisfying, for all $\r\in\Comp_p(n)$, $x\in\P(n)^{\∑_r}$, and $a_1,\dots,a_p\in A$,
    \[
      d\left(\beta_{x,\r}(a_1,\dots,a_p)\right)=\sum_{i=1}^p\beta_{x,\r\circ_i(r_i-1,1)}(a_1,\dots,a_{i},d(a_i),a_{i+1},\dots,a_p),
    \]
    where $\r\circ_i(r_i-1,1)=(r_1,\dots,r_{i-1},r_i-1,1,r_{i+1},\dots,r_p)$.
  \end{defi}
  \begin{rema}
    If $d$ is a divided $\P$-derivation of $A$, then $d$ is a $\P$-derivation for the underlying $\P$-algebra structure on $A$ given by restriction through the trace map. Indeed, recall from Remark \ref{remaPtr} that the $\P$-algebra structure is given for all $\mu\in\P(n)$, $a_1,\dots,a_n\in A$, by
    \[
      \mu(a_1,\dots,a_n)=\beta_{\mu,(1,\dots,1)}(a_1,\dots,a_n).
    \]
    From the definition of a divided $\P$-derivation, it is easy to check that $d$ is a $\P$-derivation for this action.
  \end{rema}

  We obtain the following result:
  \begin{theo}\label{theodivder}
    A divided power $\Der_\P$-algebra is a divided power $\P$-algebra endowed with a divided $\P$-derivation.
  \end{theo}
  \begin{coro}\label{expartdivder}
    A divided power $\Der_{\Com}$-algebra is a divided power algebra in the classical sense $(A,*,\gamma)$, endowed with a linear morphism $d:A\to A$ satisfying, for all $x,y\in A$, $n>0$, $d(x*y)=d(x)*y+x*d(y)$ and $d(\gamma_n(x))=\gamma_{n-1}(x)*d(x)$.
    \end{coro}
  \begin{proof}
    The beginning of this proof is analogue to the proof of Theorem \ref{theodivshift}. Following Proposition \ref{propdistrgamma}, a $\Der_\P$-algebra is a vector space $A$ endowed with the structures of a $\Gamma(\P)$-algebra and of a $\Gamma(\D)$-algebra. To understand the compatibility relations between $d$ and the $\Gamma(\P)$-algebra structure, let us compute the image of $d(\beta_{x,\r}(a_1,\dots,a_p))$ under the distributive law $\trho:\Q\tcirc(\P\tcirc V)\to \P\tcirc(\Q\tcirc V)$. 
    Applying Lemma \ref{lemmtLambda} to the element
      \[
        \mathfrak s=d(\beta_{x,\r}(a_1,\dots,a_p))=\beta_{d,(1)}(\beta_{x,\r}(a_1,\dots,a_p)),
      \]
      we obtain:
      \[
        \trho(\mathfrak s)=\as_{\P,\D,A}\circ (Tr_{\P,\D}\tcirc A)\left(\beta_{\sum_{\sigma\in\∑_{(1)\diamond \r}/\∑_{(1)}\wr\∑_{\r}} \rho(d;x)\cdot\sigma,(1)\diamond(\r)}(a_1,\dots,a_p)\right).
      \]
      Note that $(1)\diamond \r=\r$, and $\∑_{(1)\diamond \r}=\∑_{(1)}\wr\∑_{\r}=\∑_{\r}$. So,
      \begin{align*}
        \trho(\mathfrak s)&=\as_{\P,\D,A}\circ(Tr_{\P,\D}\tcirc A)\left(\beta_{\sum_{i=1}^n(x;(d^0)^{\times i-1},d,(d^0)^{\times n-i-1}),\r}(a_1,\dots,a_p)\right).\\
        &=\as_{\P,\D,A}\left(\beta_{\sum_{i=1}^n\beta_{x,(1,\dots,1)}((d^0)^{\times i-1},d,(d^0)^{\times n-i-1}),\r}(a_1,\dots,a_p)\right).\\
      \end{align*}
      According to relation \ref{relrepet} and \ref{rellin}, this is equal to
      \begin{multline*}
         \as_{\P,\D,A}\bigg(\sum_{j=1}^p\beta_{\beta_{x,(1,\dots,1)}((d^0)^{\times r_1+\dots+r_{j-1}},d,(d^0)^{\times r_j-1+r_{j+1}+\dots+r_p}),(r_1,\dots,r_{j-1},1,r_j,\dots,r_p)}\\(a_1,\dots,a_{i-1},a_i,a_i,\dots,a_p)\bigg),
      \end{multline*}
      which is equal to
      \[
        \sum_{j=1}^p\beta_{x,(r_1,\dots,r_{j-1},1,r_j,\dots,r_p)}(a_1,\dots,a_{i-1},da_i,a_i,\dots,a_p),
      \]
      which proves that $d$ is a divided $\P$-derivation of $A$.
  \end{proof}
    \begin{rema}
      In \cite{HC}, Cartan suggests that, if a divided power algebra $A$ is endowed with a differential $d$, then $d$ should satisfy the relation $d(\gamma_n(x))=\gamma_{n-1}(x)*d(x)$. Corollary \ref{expartdivder} can be seen as a confirmation of this observation.

  In \cite{KP}, Keigher and Pritchard study, for any commutative ring $A$, a divided power ring $HA$ of formal divided power series, also called `Hurwitz series', with coefficients in $A$. This ring is endowed with a derivation $\p_A:HA\to HA$.

  The augmentation ideal $H_0A$ of this ring is a divided power algebra, and the derivation $\p_A$ restricts to a divided derivation $\p_A:H_0A\to H_0A$ which turns $H_0A$ into a divided power $\Der_{\Com}$-algebra.
\end{rema}
	\section{Divided \texorpdfstring{$p$}{p}-level algebras}\label{secdivplev}
	The operad $\Lev_p$ of $p$-level algebra was introduced in \cite{SI3} as a generalisation of the operad $\Lev$ of level algebras, due to Chataur--Livernet \cite{CL}, to study unstable modules over the Steenrod algebra. The aim of this section is to characterise divided power algebras over the operad $\Lev_p$, when the base field is of characteristic $p$. This very succinct characterisation, obtained in Theorem \ref{theolevpdiv}, turns out to be very similar to the characterisation of divided power algebras in the classical sense due to Soublin \cite{JPSou} which we described in Example \ref{exdivop}. To prove this result, we use the fact that the operad $\Lev_p$ injects into the operad $\bar\Pi$ of (non-empty) partitions, and that $\bar\Pi$ is isomorphic to the operad $\Shift_{\Com}$. The divided power structure on $\Shift_{\Com}$ has been identified as an application of Theorem \ref{theodivshift} in Corollary \ref{expartdiv}.
	\begin{theo}\label{theolevpdiv}
		Let $p$ be the characteristic of the base field $\F$. A divided power $p$-level $\F$-algebra is a $p$-level algebra $(A,\star)$ endowed with an operation $\varphi_p:A\to A$ satisfying, for all $a,b\in A$, $\lambda\in\F$,
		\begin{enumerate}[label=\arabic*)]
			\item \label{rel2L1}$\varphi_p(\lambda a)=\lambda^p\varphi_p(a)$,
			\item \label{rel2L2}$\star(a^{\times p})=0$,
			\item \label{rel2L3}$\varphi_p(a+b)=\varphi_p(b)+\left(\sum_{i=1}^{p-1}\frac{\star(a^{\times i},b^{\times{p-i}})}{i!(p-i)!}\right) +\varphi_p(a)$,
			\item \label{rel2L4}$\varphi_p(\star(a_1,\dots,a_p))=0$.
		\end{enumerate}
	\end{theo}
	\begin{rema}\label{remasubopdiv}
		Let $f:\P\to \Q$ be a morphism of operads. Following Theorem \ref{theoinv}, a divided power algebra over the operad $f(\P)\subset\Q$ is a vector space endowed with the divided power operations $\beta_{f(x),\r}$ for $x\in f^{-1}(f(\P(n))^{\∑_\r})$ satisfying the divided power relations \ref{relperm} to \ref{relcomp}. In particular, if $f$ is injective, we can identify $\P$ with $f(\P)$, and since $f$ induces isomorphisms $\P(n)^{\∑_\r}\to f(\P(n))^{\∑_\r}$, a divided power $\P$-algebra is a vector space endowed with all operations $\beta_{y,\r}$ with $y\in f(\P(n)^{\∑_\r})$. 

		In other words, a divided power algebra over a suboperad of an operad $\Q$ can be defined as a vector space endowed with certain of the divided power operations of $\Q$.
	\end{rema}
	\begin{proof}[Proof of Theorem \ref{theolevpdiv}]
		Recall from Example \ref{exdivop} c) that a $\Gamma(\Com)$-algebra in characteristic $p$ is a vector space $A$ endowed with a commutative associative multiplication $*$, and with a divided $p$-th power $\gamma_p$, satisfying the relations \ref{rel2Com1} to \ref{rel2Com4}. The operations $*$ and $\gamma_p$ correspond respectively to the operations $\beta_{X_2,(1,1)}$ and $\beta_{X_p,(p)}$. From Theorem \ref{theodivshift} and the subsequent Corollary \ref{expartdiv}, we deduce that a $\Gamma(\Shift_{\Com})$-algebra in characteristic $p$ is a vector space $A$ endowed with a commutative associative multiplication $*$, a divided $p$-th power $\gamma_p$, and a linear endomorphism $d=\beta_{d,(1)}$, satisfying $*(d,d)=d\circ*$ and $\gamma_p\circ d=d\circ \gamma_p$. 

    Note that, for any vector space $V$, the vector space $\Gamma(\Shift_{\Com})(V)$ is generated by the elements of the form $\beta_{X_n,\r}(d^{i_1}v_1,\dots,d^{i_s}v_s)$ where $n\in\N$, $\r\in\Comp_s(n)$, and for all $j\in[s]$, $i_j\in\N$ and $v_j\in V$, and this is equal to:
		\[
			\prod_{j\in[s]}\prod_{k=1}^{\c(r_j)}\frac{\gamma_p^{\circ \E(r_j,k)}(d^{i_j}v_j)^{*\lambda_{r_j,k}}}{\lambda_{r_j,k}!},
		\]
		where $r_j=\sum_{k=1}^{\c(r_j)}\lambda_{r_j,k}p^{\E(r_j,k)}$ is the $p$-ary expansion of $r_j$, and the product is $*$.

		Recall from \cite{SI3} that there is an injection of operads $\iota\colon\Lev_p\to\bar\Pi$ mapping $\star$ to $(\emptyset,[p])$, whose image is spanned by the partitions $I$ satisfying the summation condition:
		\[
		 	\sum_{j\in\N}\frac{|I_j|}{p^j}=1.
		\]
		Under the identification $\bar\Pi\cong \Shift_{\Com}$, the injection $\iota$ sends $\star$ to $(X_p;d^{\times p})$, and the image of $\Lev_p$ is spanned by the $(X_n;d^{i_1},\dots,d^{i_n})$ satisfying the summation condition:
		\[
		 	\sum_{j=1}^n\frac{1}{p^{i_j}}=1.
		\]
		Following Remark \ref{remasubopdiv}, a divided power $\Lev_p$-algebra can be considered as a vector space containing those divided power operations $\beta_{x,\r}$ such that $x$ is in the image of $\Lev(p)^{\∑_\r}$ in $\Shift_{\Com}(p)^{\∑_\r}$. Let $\star$ be the operation $\beta_{(X_p;d,\dots,d);(1,\dots,1)}$, and $\varphi_p=\beta_{(X_p;d,\dots,d);(p)}$. Since $(X_p;d,\dots,d)\in\iota(\Lev_p)$, all $\Gamma(\iota(\Lev_p))$ are endowed with the operations $\star$ and $\varphi_p$. From this definition it is clear that all operations of the form
		\[
		 	u_1,\dots,u_m\mapsto \beta_{X_m,(1,\dots,1)}(d^{\alpha_1}u_1,\dots,d^{\alpha_m}u_m)
		 \]
		 with $\sum_{j=1}^m\frac{1}{p^{\alpha_j}}$ can be expressed using only $\star$. Since $\gamma_p$ and $*$ satisfy the only relations \ref{rel2Com1} to \ref{rel2Com4}, and since $*(d,d)=d\circ*$ and $\gamma_p\circ d=d\circ \gamma_p$, it is clear that $\star$ and $\varphi_p$ satisfy the only relations \ref{rel2L1} to \ref{rel2L4}. We have to prove that these two operations are enough to describe the structure of a $\Gamma(\Lev_p)$-algebra.

		From what precedes, the vector space $\Gamma(\iota(\Lev_p))(V)$ is generated by the elements of the form $\t=\beta_{X_n,\r}(d^{i_1}v_1,\dots,d^{i_s}v_s)$ for $r_1,\dots,r_s\in\N$, and $i_1,\dots,i_s\in\N$ satisfying:
		\[
		 	\sum_{j=1}^{s}\frac{r_j}{p^{i_j}}=1,
		\]
		and one has:
		\[
			\t=\prod_{j\in[s]}\prod_{k=1}^{\c(r_j)}\frac{\gamma_p^{\circ \E(r_j,k)}(d^{i_j}v_j)^{*\lambda_{r_j,k}}}{\lambda_{r_j,k}!}.
		\]
		We have to prove that these elements can be expressed using only the $p$-level multiplication $\star$, the divided $\star$-power $\varphi_p$, and the vectors of $V$. We can re-write the element $\t$:
		\[
			\t=\prod_{j\in[s]}\prod_{k=1}^{\c(r_j)}\frac{\varphi_p^{\circ \E(r_j,k)}(d^{i_j-\E(r_j,k)}v_j)^{*\lambda_{r_j,k}}}{\lambda_{r_j,k}!},
		\]
		and,
		\begin{multline*}
			\t=\beta_{X_m,(1,\dots,1)}\bigg(\frac{\varphi_p^{\circ \E(r_1,1)}(d^{i_1-\E(r_1,1)}v_1)^{*\lambda_{r_1,1}}}{\lambda_{r_1,1}!},\dots,\frac{\varphi_p^{\circ \E(r_1,\c(r_1))}(d^{i_1-\E(r_1,\c(r_1))}v_1)^{*\lambda_{r_1,\c(r_1)}}}{\lambda_{r_1,\c(r_1)}!},\\
			\dots,\frac{\varphi_p^{\circ \E(r_s,\c(r_s))}(d^{i_s-\E(r_s,\c(r_s))}v_s)^{*\lambda_{r_s,\c(r_s)}}}{\lambda_{r_s,\c(r_s)}!}\bigg),
		\end{multline*}
		where $m=\sum_{j=1}^s\c(r_j)$. We also have:
		\begin{multline*}
			\t=\left(\prod_{j\in[s],k\in\c(r_j)}\frac{1}{\lambda_{r_j,k}!}\right)\beta_{X_m,(1,\dots,1)}\bigg(\left(d^{i_1-\E(r_1,1)}\varphi_p^{\circ \E(r_1,1)}\left(v_1\right)\right)^{\times\lambda_{r_1,1}},\\
			\dots,\left(d^{i_1-\E(r_1,\c(r_1))}\varphi_p^{\circ \E(r_1,\c(r_1))}\left(v_1\right)\right)^{\times\lambda_{r_1,\c(r_1)}},
			\dots,\\\left(d^{i_s-\E(r_s,\c(r_s))}\varphi_p^{\circ \E(r_s,\c(r_s))}\left(v_s\right)\right)^{\times\lambda_{r_s,\c(r_s)}}\bigg),
		\end{multline*}
		where $m=\sum_{j\in[s],k\in\c(r_j)}\lambda_{r_j,k}$. Finally, the operation sending $(u_{j,k,l})_{j\in[s],k\in[\c(r_j)],l\in[\lambda_{r_j,k}]}$ to
		\begin{multline*}
			\beta_{X_m,(1,\dots,1)}\bigg(d^{i_1-\E(r_1,1)}(u_{1,1,1}),\dots,d^{i_1-\E(r_1,1)}(u_{1,1,\lambda_{r_1,1}}),\\\dots,d^{i_1-\E(r_1,\c(r_1))}(u_{1,\c(r_1),\lambda_{r_1,1}}),\dots, d^{i_s-\E(r_s,\c(r_s))}u_{s,\c(r_s),\lambda_{s,\c(r_s)}}\bigg)
		\end{multline*}
		can be expressed only using $\star$, because:
		\begin{align*}
			\sum_{j\in[s],k\in[\c(r_j)],l\in[\lambda_{r_j,k}]}\frac{1}{p^{i_j-\E(r_j,k)}}&=\sum_{j\in[s],k\in[\c(r_j)]}\frac{\lambda_{r_j,k}}{p^{i_j-\E(r_j,k)}}\\
			&=\sum_{j\in[s],k\in[\c(r_j)]}\frac{\lambda_{r_j,k}p^{\E(r_j,k)}}{p^{i_j}}\\
			&=\sum_{j\in[s]}\frac{r_j}{p^{i_j}}=1.
		\end{align*}
		This also proves, \textit{a posteriori}, that $i_j-\E(r_j,k)\ge0$ for all $j,k$.
	\end{proof}
  \section{Divided power Poisson algebras}\label{secPoisdiv}
  In \cite{BF}, Fresse obtained a characterisation for divided power Poisson algebras in characteristic 2, as well as a partial characterisation for divided power Poisson algebras in any prime characteristic $p>0$. The only relation which is not explicit in this result involves $\Lie$-polynomials $\Gamma_{i,j}$ that appear in the work of Cohen on the homology of $C_{n+1}$ spaces (see \cite{FrC}). The expression of these polynomials is not known in general and needs to be computed separately in each characteristic. In this section, we compute the $\Lie$-polynomials $\Gamma_{ij}$ in characteristic 3. This allows us to obtain the characterisation for divided power Poisson algebra in characteristic 3.
    \begin{theo}[\cite{BF}]\label{theodivPois}
    In characteristic $p$ prime, a $\Gamma(\Pois)$-algebra is a vector space $A$ endowed with a commutative associative multiplication $*$, a Lie bracket $[-;-]$, a divided $p$-th power operation $\gamma_p$, a Frobenius $F$, such that $(A,*,\gamma_p)$ is a divided power algebra in the classical sense, $(A,[-;-],F)$ is a $p$-restricted Lie algebra, $(A,*,[-;-])$ is a Poisson algebra, satisfying the following compatibility relations, for all $a,b\in A$:
    \begin{enumerate}[label=(P\arabic*)]
      \item\label{relP1} $[\gamma_p(a);b]=\frac{a^{* p-1}}{(p-1)!}*[a;b]$,
      \item\label{relP2} $F(a*b)=\sum_{i,j\in\N, i+j\le p}a^{*i}*b^{*j}*\Gamma_{i,j}$,
      \item\label{relP3} $F(\gamma_p(a))=0$,
    \end{enumerate}
    where the $\Gamma_{i,j}$ are products of Lie polynomials in $a$ and $b$.
  \end{theo}
  \begin{theo}[\cite{BF}]
    In characteristic $p=2$, the products of $\Lie$-polynomials $\Gamma_{i,j}$ are: $\Gamma_{0,0}=\Gamma_{1,0}=\Gamma_{0,1}=0$, and $\Gamma_{1,1}=[a;b]$, so (P2) is equivalent to $F(a*b)=a*b*[a;b]$.
  \end{theo}
\begin{theo}\label{theoPois3}
  In characteristic $p=3$, the products of $\Lie$-polynomials $\Gamma_{i,j}$ are given by: $\Gamma_{0,0}=\Gamma_{1,0}=\Gamma_{0,1}=0$, $\Gamma_{1,1}=-[a;b]^{2}$, $\Gamma_{1,2}=\left[[b;a];a\right]$, and $\Gamma_{2,1}=\left[[a;b];b\right]$. A $\Gamma(\Pois)$-algebra is a vector space $A$ endowed with a commutative associative multiplication $*$, a Lie bracket $[-;-]$, a divided cube operation $\gamma_3$, a Frobenius $F$, such that $(A,*,\gamma_3)$ is a divided power algebra in the classical sense, $(A,[-;-],F)$ is a $3$-restricted Lie algebra, $(A,*,[-;-])$ is a Poisson algebra, satisfying the compatibility relations which include the following, for all $a,b\in A$:
    \begin{enumerate}[label=(P\arabic*)]
      \item\label{rel3P1} $[\gamma_3(a);b]=\frac{a^{*2}}{2}*[a;b]$,
      \item\label{rel3P2} $F(a*b)=\left[[b;a];a\right]*a*b^{*2}+\left[[a;b];b\right]*a^{*2}*b-a*b*[a;b]^{*2}$,
      \item\label{rel3P3} $F(\gamma_3(a))=0$,
    \end{enumerate}
\end{theo}
The proof of Theorem \ref{theoPois3} is contained in Appendix \ref{appPois3}.

In conclusion, the techniques that have been developed in this article formalise the method used by Fresse to obtain the characterisation of divided power Poisson algebras, and generalise his approach to any divided power algebra structure over a product of operad with distributive law. Our new methods have been able to solve characterisation problems both in general settings, such as the case of the operad of $\P$-algebra equipped with a $\P$-endomorphism or $\P$-derivation for any operad $\P$ over a field of unspecified characteristic (see Theorems \ref{theodivshift} and \ref{theodivder}), and to obtain a more practical characterisation in a specific setting, namely the case of the operad $\Lev_p$ of $p$-level algebras over a field of characteristic $p$ (see Theorem \ref{theolevpdiv}).

In the case of the operad of Poisson algebra, this method helped us clarifying the combinatorics involved in the partial characterisation given by Fresse. However, it does not simplify the computations needed to obtain a complete characterisation. Proof of this affirmation is the somewhat lengthy computation contained in Appendix \ref{appPois3}. We do not expect that the computation needed for a similar characterisation in characteristic $p>3$, even in the case $p=5$, could be laid out in a similar way to what is done in Appendix \ref{appPois3}, without more powerful computational results. We believe, however, that techniques of computer algebra, such as those used by Bremner--Markl \cite{BrMa} and Bremner--Dotsenko, might be able to run the computation successfully in characteristic $5$ or higher.
\appendix
\section{Proof of Theorem \ref{theoPois3}}\label{appPois3}
\begin{proof}[proof of Theorem \ref{theoPois3}]
We just have to check that 
\[
  \trho(F(a*b))=\left[[b;a];a\right]*a*b^{*2}+\left[[a;b];b\right]a^{*2}*b-a*b*[a;b]^{*2}.
\]
Recall that $F(a*b)=\beta_{L_3,(3)}(\beta_{X_2,(1,1)}(a,b))$, where $L_3=[-,-]\circ_1[-;-]\cdot(\id_{[3]}+(2\ 3))$. According to Lemma \ref{lemmtLambda}, one has:
  \begin{multline*}
    \trho\left(\beta_{L_3,(3)}(\beta_{X_2,(1,1)}(a,b))\right)=\\\as_{\Com,\Lie,A}\circ(Tr_{\Com,\Lie}\tcirc A)\left(\beta_{\sum_{\sigma\in\∑_{(3)\diamond(1,1)}/\∑_{3}\wr\∑_{(1,1)} }\rho(L_3;X_2^{\times 3})\cdot \sigma,(3)\diamond\left(1,1\right)}\left(a,b\right)\right).
  \end{multline*}
From now on, the multiplication $*$ will be denoted by a juxtaposition of letters, $a*b=ab$.

We will compute $\rho(L_3,X_2^{\times 3})$ using formal variables: if $\mu\in \P(n)$ and $\sigma\in\∑_n$, we can denote by $\mu(a_{\sigma(1)},\dots,a_{\sigma(n)})=\mu\cdot\sigma$. The distributive law $\rho$ can be denoted by rewriting: $[a_1a_2;a_3]\leadsto a_1[a_2;a_3]+a_2[a_1;a_3]$. Using $[a_1;a_2]=-[a_2;a_1]$, we also have $[a_1;a_2a_3]\leadsto[a_1;a_2]a_3+[a_1;a_3]a_2$. We get:
  \begin{align*}
   \hskip-20pt \left[[a_1a_2;a_3a_4];a_5a_6\right]&\leadsto [a_1[a_2;a_3a_4];a_5a_6]+[a_2[a_1;a_3a_4]a_5a_6],\\
    &\leadsto a_1[[a_2;a_3a_4];a_5a_6]+[a_2;a_3a_4][a_1;a_5a_6]\\
    &+a_2[[a_1;a_3a_4];a_5a_6]+[a_1;a_3a_4][a_2;a_5a_6],\\
    &\leadsto a_1[[a_2;a_3a_4];a_5]a_6+ a_1[[a_2;a_3a_4];a_6]a_5+[a_2;a_3a_4][a_1;a_5]a_6+[a_2;a_3a_4][a_1;a_6]a_5\\
    &+a_2[[a_1;a_3a_4];a_5]a_6+a_2[[a_1;a_3a_4];a_6]a_5+[a_1;a_3a_4][a_2;a_5]a_6+[a_1;a_3a_4][a_2;a_6]a_5,\\
    &\leadsto a_1[[a_2;a_3]a_4;a_5]a_6+a_1[[a_2;a_4]a_3;a_5]a_6+a_1[[a_2;a_3]a_4;a_6]a_5+a_1[[a_2;a_4]a_3;a_6]a_5\\
    &+[a_2;a_3]a_4[a_1;a_5]a_6+[a_2;a_4]a_3[a_1;a_5]a_6+[a_2;a_3]a_4[a_1;a_6]a_5+[a_2;a_4]a_3[a_1;a_6]a_5\\
    &+a_2[[a_1;a_3]a_4;a_5]a_6+a_2[[a_1;a_4]a_3;a_5]a_6+a_2[[a_1;a_3]a_4;a_6]a_5+a_2[[a_1;a_4]a_3;a_6]a_5\\
    &+[a_1;a_3]a_4[a_2;a_5]a_6+[a_1;a_4]a_3[a_2;a_5]a_6+[a_1;a_3]a_4[a_2;a_6]a_5+[a_1;a_4]a_3[a_2;a_6]a_5,\\
    &\leadsto a_1[a_2;a_3][a_4;a_5]a_6+a_1[[a_2;a_3];a_5]a_4a_6+a_1[a_2;a_4][a_3;a_5]a_6+a_1a_3[[a_2;a_4];a_5]a_6\\
    &+a_1[a_2;a_3][a_4;a_6]a_5+a_1a_4[[a_2;a_3];a_6]a_5+a_1[a_2;a_4][a_3;a_6]a_5+a_1a_3[[a_2;a_4];a_6]a_5\\
    &+[a_2;a_3]a_4[a_1;a_5]a_6+[a_2;a_4]a_3[a_1;a_5]a_6+[a_2;a_3]a_4[a_1;a_6]a_5+[a_2;a_4]a_3[a_1;a_6]a_5\\
    &+a_2[a_1;a_3][a_4;a_5]a_6+a_2a_4[[a_1;a_3];a_5]a_6+a_2[a_1;a_4][a_3;a_5]a_6+a_2a_3[[a_1;a_4];a_5]a_6\\
    &+a_2[a_1;a_3][a_4;a_6]a_5+a_2a_4[[a_1;a_3];a_6]a_5+a_2[a_1;a_4][a_3;a_6]a_5+a_2a_3[[a_1;a_4];a_6]a_5\\
    &+[a_1;a_3]a_4[a_2;a_5]a_6+[a_1;a_4]a_3[a_2;a_5]a_6+[a_1;a_3]a_4[a_2;a_6]a_5+[a_1;a_4]a_3[a_2;a_6]a_5.\\
  \end{align*}
  Then, since $(L_3;X_2^{\times 3})=([-;-]\circ_1[-;-];X_2^{\times 3})\cdot(\id_{[6]}+(3\ 5)(4 6))$, we can rewrite $(L_3;X_2^{\times 3})$ as
  \begin{align*}
    &a_1[a_2;a_3][a_4;a_5]a_6-a_1[a_2;a_5][a_3;a_6]a_4+a_1L_3(a_2,a_3,a_5)a_4a_6+a_1[a_2;a_4][a_3;a_5]a_6\\
    &-a_1[a_2;a_6][a_3;a_5]a_4+a_1a_3[[a_2;a_4];a_5]a_6+a_1a_5[[a_2;a_6];a_3]a_4+a_1[a_2;a_3][a_4;a_6]a_5\\
    &-a_1[a_2;a_5][a_4;a_6]a_3+a_1a_4[[a_2;a_3];a_6]a_5+a_1a_6[[a_2;a_5];a_4]a_3+a_1[a_2;a_4][a_3;a_6]a_5\\
    &-a_1[a_2;a_6][a_4;a_5]a_3+a_1a_3L_3(a_2,a_4,a_6)a_5+[a_2;a_3]a_4[a_1;a_5]a_6+[a_2;a_5]a_4[a_1;a_3]a_6\\
    &+[a_2;a_4]a_3[a_1;a_5]a_6+[a_2;a_6][a_1;a_3]a_4a_5+[a_2;a_3]a_4[a_1;a_6]a_5+[a_2;a_5]a_3[a_1;a_4]a_6\\
  &+[a_2;a_4]a_3[a_1;a_6]a_5+[a_2;a_6]a_5[a_1;a_4]a_3+a_2[a_1;a_3][a_4;a_5]a_6-a_2[a_1;a_5][a_3;a_6]a_4\\
    &+a_2a_4L_3(a_1,a_3,a_5)a_6+a_2[a_1;a_4][a_3;a_5]a_6-a_2[a_1;a_6][a_3;a_5]a_4+a_2a_3[[a_1;a_4];a_5]a_6\\
    &+a_2a_5[[a_1;a_6];a_3]a_6+a_2[a_1;a_3][a_4;a_6]a_5-a_2[a_1;a_5][a_4;a_6]a_3+a_2a_4[[a_1;a_3];a_6]a_5\\
    &+a_2a_6[[a_1;a_5];a_4]a_3+a_2[a_1;a_4][a_3;a_6]a_5-a_2[a_1;a_6][a_4;a_5]a_3+a_2a_3L_3(a_1,a_4,a_6)a_5\\
    &+[a_1;a_3]a_4[a_2;a_5]a_6+[a_1;a_5]a_4[a_2;a_3]a_6+[a_1;a_4]a_3[a_2;a_5]a_6+[a_1;a_6]a_5[a_2;a_3]a_4\\
    &+[a_1;a_3]a_4[a_2;a_6]a_5+[a_1;a_5]a_6[a_2;a_4]a_3+[a_1;a_4]a_3[a_2;a_6]a_5+[a_1;a_6]a_5[a_2;a_4]a_3,
  \end{align*}
  which is equal to
    \begin{align*}
     \tag{$l_1$}\label{L1}&a_1\left(L_3(a_2,a_3,a_5)a_4a_6+L_3(a_2,a_4,a_5)a_3a_6+L_3(a_2,a_3,a_6)a_4a_5+L_3(a_2,a_4,a_6)a_3a_5\right)\\
     \notag\\
     \notag&+a_2\Big(L_3(a_1,a_3,a_5)a_2a_4a_6+L_3(a_1,a_4,a_5)a_2a_3a_6\\
      \tag{$l_2$}\label{L2}&+L_3(a_1,a_3,a_6)a_2a_4a_5+L_3(a_1,a_4,a_6)a_2a_3a_5\Big)\\
     \notag\\
     \notag&+[a_3;a_2][a_5;a_4]a_1a_6+[a_5;a_2][a_3;a_6]a_1a_4+[a_1;a_6][a_3;a_2]a_4a_5\\
     \tag{$l_3$}\label{L3}&+[a_1;a_4][a_5;a_2]a_3a_6+[a_1;a_4][a_3;a_6]a_2a_5+[a_1;a_6][a_5;a_4]a_2a_3\\
    \notag\\
    \notag&+[a_3;a_5][a_2;a_4]a_1a_6+[a_3;a_5][a_6;a_2]a_1a_4+[a_1;a_5][a_4;a_2]a_3a_6\\
    \tag{$l_4$}\label{L4}&+[a_1;a_3][a_6;a_2]a_4a_5+[a_1;a_3][a_4;a_6]a_2a_5+[a_1;a_5][a_6;a_4]a_2a_3\\
    \notag\\
    &\notag+[a_2;a_3][a_4;a_6]a_1a_5-[a_2;a_5][a_4;a_6]a_1a_3-[a_6;a_3][a_2;a_4]a_1a_5\\
    &\tag{$l_5$}\label{L5}-[a_4;a_5][a_2;a_6]a_1a_3+[a_6;a_1][a_2;a_4]a_3a_5+[a_4;a_1][a_2;a_6]a_3a_5\\
    \notag\\
    &\notag-[a_1;a_5][a_2;a_3]a_4a_6-[a_1;a_3][a_2;a_5]a_4a_6+[a_1;a_3][a_4;a_5]a_2a_6\\
    &\tag{$l_6$}\label{L6}+[a_1;a_5][a_6;a_3]a_2a_4-[a_3;a_5][a_4;a_1]a_2a_6+[a_3;a_5][a_6;a_1]a_2a_4.
  \end{align*}
  We will denote $(L_3;X_2^{\times 3})=\sum_{i=1}^6l_i$ where $l_i$ are the groups of summands as above.

  Observe that $(3)\diamond(1,1)$ is the partition $(\{1,3,5\},\{2,4,6\})$, and $\∑_3\wr\∑_{(1,1)}$ is the following subgroup of $\∑_6$:
  \[
    \left\{\id_{[6]}, (1\ 3)(2\ 4), (1\ 5)(2\ 6), (3\ 5)(4\ 6), (1\ 3\ 5)(2\ 4\ 6), (1\ 5\ 3)(2\ 6\ 4)\right\},
  \]
  where $\id_{[n]}\in\∑_n$ denotes the neutral element, and we denoted the other permutations as a product of disjoint cycles. Then $E=\{id_{[6]},(1\ 3), (3\ 5), (1\ 3)(2\ 4\ 6), (1\ 3\ 5), (2\ 4\ 6)\}$ is a set of representative for $\∑_{(3)\diamond(1,1)}/\∑_3\wr\∑_{(1,1)}$.

  Let us compute $\sum_{\sigma\in E}l_i\sigma$, for $i\in[6]$:
  \begin{align*}
    \sum_{\sigma\in E}l_1\sigma&=a_1\left(L_3(a_2,a_3,a_5)a_4a_6+L_3(a_2,a_4,a_5)a_3a_6+L_3(a_2,a_3,a_6)a_4a_5+L_3(a_2,a_4,a_6)a_3a_5\right)\\
    &+a_3\left(L_3(a_2,a_1,a_5)a_4a_6+L_3(a_2,a_4,a_5)a_1a_6+L_3(a_2,a_1,a_6)a_4a_5+L_3(a_2,a_4,a_6)a_1a_5\right)\\
    &+a_1\left(L_3(a_2,a_5,a_3)a_4a_6+L_3(a_2,a_4,a_3)a_5a_6+L_3(a_2,a_5,a_6)a_4a_3+L_3(a_2,a_4,a_6)a_3a_5\right)\\
    &+a_3\left(L_3(a_4,a_1,a_5)a_6a_2+L_3(a_4,a_6,a_5)a_1a_2+L_3(a_4,a_1,a_2)a_6a_5+L_3(a_4,a_6,a_2)a_1a_5\right)\\
    &+a_3\left(L_3(a_2,a_5,a_1)a_4a_6+L_3(a_2,a_4,a_1)a_5a_6+L_3(a_2,a_5,a_6)a_4a_1+L_3(a_2,a_4,a_6)a_5a_1\right)\\
    &+a_1\left(L_3(a_4,a_3,a_5)a_6a_2+L_3(a_4,a_6,a_5)a_3a_2+L_3(a_4,a_3,a_2)a_6a_5+L_3(a_4,a_6,a_2)a_3a_5\right)\\
  \end{align*}
  Since $L_3$ is fixed by all permutations in $\∑_3$, this implies that:
  \begin{align*}
    \sum_{\sigma\in E}l_1\sigma&=a_1\Big(L_3(a_2,a_3,a_6)a_4a_5+L_3(a_3,a_4,a_5)a_2a_6-L_3(a_2,a_3,a_5)a_4a_6-L_3(a_2,a_4,a_5)a_3a_6\\
    &-L_3(a_2,a_3,a_4)a_5a_6-L_3(a_2,a_5,a_6)a_3a_4-L_3a(a_4,a_5,a_6)a_2a_3\Big)\\
    &+a_3\left(L_3(a_1,a_2,a_6)a_4a_5+L_3(a_1,a_4,a_5)a_3a_6-L_3(a_1,a_2,a_4)a_5a_6-L_3(a_1,a_2,a_5)a_4a_6\right).
  \end{align*}

  \begin{align*}
  \sum_{\sigma\in E}l_2\sigma&=a_2\left(L_3(a_1,a_3,a_5)a_4a_6+L_3(a_1,a_4,a_5)a_3a_6+L_3(a_1,a_3,a_6)a_4a_5+L_3(a_1,a_4,a_6)a_3a_5\right)\\
  &+a_2\left(L_3(a_3,a_1,a_5)a_4a_6+L_3(a_3,a_4,a_5)a_1a_6+L_3(a_3,a_1,a_6)a_4a_5+L_3(a_3,a_4,a_6)a_1a_5\right)\\
  &+a_2\left(L_3(a_1,a_5,a_3)a_4a_6+L_3(a_1,a_4,a_3)a_5a_6+L_3(a_1,a_5,a_6)a_4a_3+L_3(a_1,a_4,a_6)a_5a_3\right)\\
  &+a_4\left(L_3(a_3,a_1,a_5)a_6a_2+L_3(a_3,a_6,a_5)a_1a_2+L_3(a_3,a_1,a_2)a_6a_5+L_3(a_3,a_6,a_2)a_1a_5\right)\\
  &+a_2\left(L_3(a_3,a_5,a_1)a_4a_6+L_3(a_3,a_4,a_1)a_5a_6+L_3(a_3,a_5,a_6)a_4a_1+L_3(a_3,a_4,a_6)a_5a_1\right)\\
  &+a_4\left(L_3(a_1,a_3,a_5)a_6a_2+L_3(a_1,a_6,a_5)a_3a_2+L_3(a_1,a_3,a_2)a_6a_5+L_3(a_1,a_6,a_2)a_3a_5\right),\\
  &=a_2\Big(L_3(a_1,a_4,a_5)a_3a_6-L_3(a_1,a_3,a_6)a_4a_5-L_3(a_1,a_4,a_6)a_3a_5-L_3(a_1,a_3,a_4)a_5a_6\\
  &-L_3(a_1,a_5,a_6)a_3a_4-L_3(a_3,a_4,a_6)a_1a_5-L_3(a_3,a_5,a_6)a_1a_4\Big)\\
  &+a_4a_5\left(L_3(a_1,a_2,a_6)a_3+L_3(a_2,a_3,a_6)a_1-L_3(a_1,a_2,a_3)a_6\right)\\
  &+a_1a_2a_6L_3(a_3,a_4,a_5).
  \end{align*}
  So,
 \begin{equation*}
    \sum_{\sigma\in E}(l_1+l_2)\sigma=-\sum_{(i_1,\dots,i_6)}L_3(a_{i_1},a_{i_2},a_{i_3})a_{i_4}a_{i_5}a_{i_6},
 \end{equation*}
 where $(i_1,\dots,i_6)$ ranges over the following set:
 \begin{multline*}
  \Big\{(1,2,3,4,5,6), (1,2,4,3,5,6), (1,2,5,3,4,6), (1,2,6,3,4,5), (1,3,4,2,5,6),(1,3,6,2,4,5),\\
    (1,4,5,2,3,6), (1,4,6,2,3,5), (1,5,6,2,3,4), (2,3,4,1,5,6), (2,3,5,1,4,6),(2,3,6,1,4,5),\\
     (2,4,5,1,3,6), (2,5,6,1,3,4), (3,4,5,1,2,6), (3,4,6,1,2,5), (3,5,6,1,2,4), (4,5,6,1,2,3)\Big\}
 \end{multline*}
  We now rewrite $\sum_{\sigma\in E}(l_1+l_2)\sigma$ as a sum of elements of $\Pois(6)^{\∑_{\{1,3,5\}}\times\∑_{\{2,4,6\}}}$:
  \[
    \sum_{\sigma\in E}(l_1+l_2)\sigma=-\sum_{(i_1,\dots,i_6)\in A_1}L_3(a_{i_1},a_{i_2},a_{i_3})a_{i_4}a_{i_5}a_{i_6}\\
    -\sum_{(i_1,\dots,i_6)\in A_2}L_3(a_{i_1},a_{i_2},a_{i_3})a_{i_4}a_{i_5}a_{i_6},
  \]
  where
  \begin{multline*}
    A_1=\Big\{(1,2,3,4,5,6), (1,2,5,3,4,6), (1,3,4,2,5,6), (1,3,6,2,4,5),\\
     (1,4,5,2,3,6),(1,5,6,2,3,4), (2,3,5,1,4,6), (3,4,5,1,2,6), (3,5,6,1,2,4)\Big\},
  \end{multline*}
  and
  \begin{multline*}
    A_2=\Big\{(1,2,4,3,5,6), (1,2,6,3,4,5), (1,4,6,2,3,5), (2,3,4,1,5,6),\\
     (2,3,6,1,4,5), (2,4,5,1,3,6), (2,5,6,1,3,4), (3,4,6,1,2,5), (4,5,6,1,2,3)\Big\}.
  \end{multline*}

  Now, observe that, since $x_4\circ_1L_3$ is fixed under the action of $\∑_3\times\∑_3$,
  \[
    \sum_{\sigma\in\∑_{\{1,3,5\}}\times \∑_{\{2,4,6\}}}\left(L_3(a_1,a_2,a_3)a_4a_5a_6\right)\sigma=4\sum_{(i_1,\dots,i_6)\in A_1}L_3(a_{i_1},a_{i_2},a_{i_3})a_{i_4}a_{i_5}a_{i_6},
  \]
  and similarly,
  \[
    \sum_{\sigma\in\∑_{\{1,3,5\}}\times \∑_{\{2,4,6\}}}\left(L_3(a_1,a_2,a_4)a_3a_5a_6\right)\sigma=4\sum_{(i_1,\dots,i_6)\in A_2}L_3(a_{i_1},a_{i_2},a_{i_3})a_{i_4}a_{i_5}a_{i_6}.
  \]
  So,
  \begin{multline*}
      \beta_{\sum_{\sigma\in\∑_{(p)\diamond(1,1)}/\∑_{p}\wr\∑_{(1,1)} }(l_1+l_2)\cdot \sigma,(p)\diamond\left(1,1\right)}\left(a,b\right)\\=-\beta_{\sum_{\sigma\in\∑_{\{1,3,5\}}\times \∑_{\{2,4,6\}}}\left(L_3(a_1,a_2,a_3)a_4a_5a_6\right)\sigma,(\{1,3,5\},\{2,4,6\})}(a,b)\\-\beta_{\sum_{\sigma\in\∑_{\{1,3,5\}}\times \∑_{\{2,4,6\}}}\left(L_3(a_1,a_2,a_4)a_3a_5a_6\right)\sigma,(\{1,3,5\},\{2,4,6\})}(a,b),
  \end{multline*}
  which, according to \ref{relrepet}, is equal to:
  \begin{equation*}
    -\beta_{L_3(a_1,a_2,a_3)a_4a_5a_6,(1,\dots,1)}(a,b,a,b,a,b)-\beta_{L_3(a_1,a_2,a_4)a_3a_5a_6,(1,\dots,1)}(a,b,a,b,a,b,a).
  \end{equation*}
  So, 
  \begin{multline*}
        \as_{\Com,\Lie,A}\circ(Tr_{\Com,\Lie}\tcirc A)\left(\beta_{\sum_{\sigma\in\∑_{(p)\diamond(1,1)}/\∑_{p}\wr\∑_{(1,1)} }(l_1+l_2)\cdot \sigma,(p)\diamond\left(1,1\right)}\left(a,b\right)\right)=\\
        \as_{\Com,\Lie,A}\bigg(-\beta_{\beta_{X_4,(1,1,1,1)}(L_3,1_{\Lie}^{\times 3}),(1,\dots,1)}(a,b,a,b,a,b)\\-\beta_{\beta_{X_4,(1,1,1,1)}(L_3,1_{\Lie}^{\times 3})(3\ 4),(1,\dots,1)}(a,b,a,b,a,b)\bigg),
  \end{multline*}
  which is equal to:
    \[
        -\beta_{X_4,(1,\dots,1)}(\beta_{L_3,(1,1,1)}(a,b,a),b,a,b)-\beta_{X_4,(1,\dots,1)}(\beta_{L_3,(1,1,1)}(a,b,b),a,a,b).
    \]   
Note that $L_3=[-;-]\left((id_{[3]}+(2\ 3))\right)$, so finally,
\begin{multline*}
        \as_{\Com,\Lie,A}\circ(Tr_{\Com,\Lie}\tcirc A)\left(\beta_{\sum_{\sigma\in\∑_{(p)\diamond(1,1)}/\∑_{p}\wr\∑_{(1,1)} }(l_1+l_2)\cdot \sigma,(p)\diamond\left(1,1\right)}\left(a,b\right)\right)=\\
        -\left[[a;b];a\right]ab^2-\left[[a;a];b\right]ab^2+\left[[a;b];b\right]a^2b,
\end{multline*}
  and with Jacobi's identity,
  \begin{multline*}
        \as_{\Com,\Lie,A}\circ(Tr_{\Com,\Lie}\tcirc A)\left(\beta_{\sum_{\sigma\in\∑_{(p)\diamond(1,1)}/\∑_{p}\wr\∑_{(1,1)} }(l_1+l_2)\cdot \sigma,(p)\diamond\left(1,1\right)}\left(a,b\right)\right)=\\
        \left[[b;a];a\right]ab^2+\left[[a;b];b\right]a^2b.
\end{multline*}
We now need to compute $\sum_{\sigma\in E}l_i\sigma$ for $i\in\{3,\dots,6\}$. 
One has:
\begin{align*}
  \sum_{\sigma\in E}l_3\sigma&=[a_3;a_2][a_5;a_4]a_1a_6+[a_5;a_2][a_3;a_6]a_1a_4+[a_1;a_6][a_3;a_2]a_4a_5\\
  &+[a_1;a_4][a_5;a_2]a_3a_6+[a_1;a_4][a_3;a_6]a_2a_5+[a_1;a_6][a_5;a_4]a_2a_3\\
  &+[a_1;a_2][a_5;a_4]a_3a_6+[a_5;a_2][a_1;a_6]a_3a_4+[a_3;a_6][a_1;a_2]a_4a_5\\
  &+[a_3;a_4][a_5;a_2]a_1a_6+[a_3;a_4][a_1;a_6]a_2a_5+[a_3;a_6][a_5;a_4]a_2a_1\\
  &+[a_5;a_2][a_3;a_4]a_1a_6+[a_3;a_2][a_5;a_6]a_1a_4+[a_1;a_6][a_5;a_2]a_4a_3\\
  &+[a_1;a_4][a_3;a_2]a_5a_6+[a_1;a_4][a_5;a_6]a_2a_3+[a_1;a_6][a_3;a_4]a_2a_5\\
  &+[a_1;a_6][a_5;a_4]a_3a_2+[a_1;a_6][a_3;a_2]a_4a_5+[a_1;a_2][a_5;a_6]a_4a_3\\
  &+[a_1;a_4][a_5;a_6]a_3a_2+[a_5;a_4][a_1;a_2]a_3a_6+[a_3;a_2][a_1;a_4]a_6a_5\\
  &+[a_5;a_2][a_1;a_4]a_3a_6+[a_1;a_2][a_5;a_6]a_3a_4+[a_3;a_6][a_5;a_2]a_4a_1\\
  &+[a_3;a_4][a_1;a_2]a_5a_6+[a_3;a_4][a_5;a_6]a_2a_1+[a_3;a_6][a_1;a_4]a_2a_5\\
  &+[a_3;a_6][a_5;a_4]a_1a_2+[a_3;a_6][a_1;a_2]a_4a_5+[a_3;a_2][a_5;a_6]a_4a_1\\
  &+[a_3;a_4][a_5;a_6]a_1a_2+[a_5;a_4][a_3;a_2]a_1a_6+[a_1;a_2][a_3;a_4]a_6a_5\\
  &=-\sum_{\sigma\in\∑_{\{1,3,5\}\times\∑_{\{2,4,6\}}}}[a_1;a_2][a_3;a_4]a_5a_6,
\end{align*}
so
\begin{multline*}
      \beta_{\sum_{\sigma\in\∑_{(p)\diamond(1,1)}/\∑_{p}\wr\∑_{(1,1)} }l_3\cdot \sigma,(p)\diamond\left(1,1\right)}\left(a,b\right)\\=-\beta_{\sum_{\sigma\in\∑_{\{1,3,5\}}\times \∑_{\{2,4,6\}}}\left([a_1;a_2][a_3;a_4]a_5a_6\right)\sigma,(\{1,3,5\},\{2,4,6\})}(a,b),
  \end{multline*}
  which, according to \ref{relrepet}, is equal to:
  
  \[
    -\beta_{[a_1;a_2][a_3;a_4]a_5a_6,(1,\dots,1)}(a,b,a,b,a,b)=-ab[a;b]^2.
  \]
  One has:
\begin{align*}
  \sum_{\sigma\in E}l_4\sigma=&[a_3;a_5][a_2;a_4]a_1a_6+[a_3;a_5][a_6;a_2]a_1a_4+[a_1;a_5][a_4;a_2]a_3a_6\\
  &+[a_1;a_3][a_6;a_2]a_4a_5+[a_1;a_3][a_4;a_6]a_2a_5+[a_1;a_5][a_6;a_4]a_2a_3\\
  &+[a_1;a_5][a_2;a_4]a_3a_6+[a_1;a_5][a_6;a_2]a_3a_4+[a_3;a_5][a_4;a_2]a_1a_6\\
  &-[a_1;a_3][a_6;a_2]a_4a_5-[a_1;a_3][a_4;a_6]a_2a_5+[a_3;a_5][a_6;a_4]a_2a_1\\
  &-[a_3;a_5][a_2;a_4]a_1a_6-[a_3;a_5][a_6;a_2]a_1a_4+[a_1;a_3][a_4;a_2]a_5a_6\\
  &+[a_1;a_5][a_6;a_2]a_4a_3+[a_1;a_5][a_4;a_6]a_2a_3+[a_1;a_3][a_6;a_4]a_2a_5\\
  &+[a_1;a_5][a_4;a_6]a_3a_2+[a_1;a_5][a_2;a_4]a_3a_6+[a_3;a_5][a_6;a_4]a_1a_2\\
  &-[a_1;a_3][a_2;a_4]a_6a_5-[a_1;a_3][a_6;a_2]a_4a_5+[a_3;a_5][a_2;a_6]a_4a_1\\
  &+[a_5;a_1][a_2;a_4]a_3a_6+[a_5;a_1][a_6;a_2]a_3a_4+[a_3;a_1][a_4;a_2]a_5a_6\\
  &+[a_3;a_5][a_6;a_2]a_4a_1+[a_3;a_5][a_4;a_6]a_2a_1+[a_3;a_1][a_6;a_4]a_2a_5\\
  &+[a_3;a_5][a_4;a_6]a_1a_2+[a_3;a_5][a_2;a_4]a_1a_6+[a_1;a_5][a_6;a_4]a_3a_2\\
  &+[a_1;a_3][a_2;a_4]a_6a_5+[a_1;a_3][a_6;a_2]a_4a_5+[a_1;a_5][a_2;a_6]a_4a_3,\\
  &=0,
\end{align*}
so
\begin{equation*}
      \beta_{\sum_{\sigma\in\∑_{(p)\diamond(1,1)}/\∑_{p}\wr\∑_{(1,1)} }l_4\cdot \sigma,(p)\diamond\left(1,1\right)}\left(a,b\right)=-\beta_{0,(\{1,3,5\},\{2,4,6\})}(a,b)=0,
  \end{equation*}
  according to \ref{rellin}.
  One has:
\begin{align*}
  \sum_{\sigma\in E}l_5\sigma&=[a_2;a_3][a_4;a_6]a_1a_5-[a_2;a_5][a_4;a_6]a_1a_3-[a_6;a_3][a_2;a_4]a_1a_5\\
  &-[a_4;a_5][a_2;a_6]a_1a_3+[a_6;a_1][a_2;a_4]a_3a_5+[a_4;a_1][a_2;a_6]a_3a_5\\
  &+[a_2;a_1][a_4;a_6]a_3a_5-[a_2;a_5][a_4;a_6]a_1a_3-[a_6;a_1][a_2;a_4]a_3a_5\\
  &-[a_4;a_5][a_2;a_6]a_1a_3+[a_6;a_3][a_2;a_4]a_1a_5+[a_4;a_3][a_2;a_6]a_1a_5\\
  &+[a_2;a_5][a_4;a_6]a_1a_3-[a_2;a_3][a_4;a_6]a_1a_5-[a_6;a_5][a_2;a_4]a_1a_3\\
  &-[a_4;a_3][a_2;a_6]a_1a_5+[a_6;a_1][a_2;a_4]a_3a_5+[a_4;a_1][a_2;a_6]a_3a_5\\
  &+[a_4;a_1][a_6;a_2]a_3a_5-[a_4;a_5][a_6;a_2]a_1a_3-[a_2;a_1][a_4;a_6]a_3a_5\\
  &-[a_6;a_5][a_4;a_2]a_1a_3+[a_2;a_3][a_4;a_6]a_1a_5+[a_6;a_3][a_4;a_2]a_1a_5\\
  &+[a_2;a_5][a_4;a_6]a_3a_1-[a_2;a_1][a_4;a_6]a_3a_5-[a_6;a_5][a_2;a_4]a_3a_1\\
  &-[a_4;a_1][a_2;a_6]a_3a_5+[a_6;a_3][a_2;a_4]a_5a_1+[a_4;a_3][a_2;a_6]a_5a_1\\
  &+[a_4;a_3][a_6;a_2]a_1a_5-[a_4;a_5][a_6;a_2]a_1a_3-[a_2;a_3][a_4;a_6]a_1a_5\\
  &-[a_6;a_5][a_4;a_2]a_1a_3+[a_2;a_1][a_4;a_6]a_3a_5+[a_6;a_1][a_4;a_2]a_3a_5,\\
  &=0,
\end{align*}
so
\begin{equation*}
      \beta_{\sum_{\sigma\in\∑_{(p)\diamond(1,1)}/\∑_{p}\wr\∑_{(1,1)} }l_5\cdot \sigma,(p)\diamond\left(1,1\right)}\left(a,b\right)=-\beta_{0,(\{1,3,5\},\{2,4,6\})}(a,b)=0,
  \end{equation*}
  according to \ref{rellin}.
  One has:
\begin{align*}
  \sum_{\sigma\in E}l_6\sigma&=-[a_1;a_5][a_2;a_3]a_4a_6-[a_1;a_3][a_2;a_5]a_4a_6+[a_1;a_3][a_4;a_5]a_2a_6\\
  &+[a_1;a_5][a_6;a_3]a_2a_4-[a_3;a_5][a_4;a_1]a_2a_6+[a_3;a_5][a_6;a_1]a_2a_4\\
  &-[a_3;a_5][a_2;a_1]a_4a_6+[a_1;a_3][a_2;a_5]a_4a_6-[a_1;a_3][a_4;a_5]a_2a_6\\
  &+[a_3;a_5][a_6;a_1]a_2a_4-[a_1;a_5][a_4;a_3]a_2a_6+[a_1;a_5][a_6;a_3]a_2a_4\\
  &-[a_1;a_3][a_2;a_5]a_4a_6-[a_1;a_5][a_2;a_3]a_4a_6+[a_1;a_5][a_4;a_3]a_2a_6\\
  &+[a_1;a_3][a_6;a_5]a_2a_4+[a_3;a_5][a_4;a_1]a_2a_6-[a_3;a_5][a_6;a_1]a_2a_4\\
  &-[a_3;a_5][a_4;a_1]a_6a_2+[a_1;a_3][a_4;a_5]a_6a_2-[a_1;a_3][a_6;a_5]a_4a_2\\
  &+[a_3;a_5][a_2;a_1]a_4a_6-[a_1;a_5][a_6;a_3]a_4a_2+[a_1;a_5][a_2;a_3]a_4a_6\\
  &-[a_3;a_1][a_2;a_5]a_4a_6-[a_3;a_5][a_2;a_1]a_4a_6+[a_3;a_5][a_4;a_1]a_2a_6\\
  &+[a_3;a_1][a_6;a_5]a_2a_4-[a_5;a_1][a_4;a_3]a_2a_6+[a_5;a_1][a_6;a_3]a_2a_4\\
  &-[a_1;a_5][a_4;a_3]a_6a_2-[a_1;a_3][a_4;a_5]a_6a_2+[a_1;a_3][a_6;a_5]a_4a_2\\
  &+[a_1;a_5][a_2;a_3]a_4a_6-[a_3;a_5][a_6;a_1]a_4a_2+[a_3;a_5][a_2;a_1]a_4a_6,\\
  &=0,
\end{align*}
so
\begin{equation*}
      \beta_{\sum_{\sigma\in\∑_{(p)\diamond(1,1)}/\∑_{p}\wr\∑_{(1,1)} }l_6\cdot \sigma,(p)\diamond\left(1,1\right)}\left(a,b\right)=-\beta_{0,(\{1,3,5\},\{2,4,6\})}(a,b)=0,
  \end{equation*}
  according to \ref{rellin}.

  Finally,
\begin{align*}
  \trho(F(a*b))&=\as_{\Com,\Lie,A}\circ(Tr_{\Com,\Lie}\tcirc A)\left(\beta_{\sum_{\sigma\in\∑_{(p)\diamond(1,1)}/\∑_{p}\wr\∑_{(1,1)} }\sum_{i=1}^6l_i\cdot \sigma,(p)\diamond\left(1,1\right)}\left(a,b\right)\right)\\
  &=\left[[b;a];a\right]*a*b^{*2}+\left[[a;b];b\right]*a^{*2}*b-a*b*[a;b]^{*2}.
\end{align*}
\end{proof}

  \vfill
Sacha {\scshape Ikonicoff}, {\scshape pims--cnrs} postdoctoral associate, University of Calgary\\
e-mail: \texttt{sacha.ikonicoff@ucalgary.ca}
\end{document}